\documentclass[default]{sn-jnl}

\usepackage{graphicx}\usepackage{multirow}\usepackage{amsmath}\usepackage{amssymb}\usepackage{amsfonts}\usepackage{amsthm}\usepackage{mathrsfs}\usepackage[title]{appendix}\usepackage{xcolor}\usepackage{textcomp}\usepackage{manyfoot}\usepackage{booktabs}\usepackage{algorithm}\usepackage{algorithmicx}\usepackage{algpseudocode}\usepackage{listings}

\usepackage{lineno}
\usepackage{hyperref}
\modulolinenumbers[5]

\RequirePackage{amsmath}
\RequirePackage{amsfonts}
\RequirePackage{bbm}
\RequirePackage[normalem]{ulem}

\newcommand{\vecs}[1]{\mathbf{#1}}

\newcommand{\comdom}{V}

\usepackage[capitalise]{cleveref}
\crefname{equation}{}{}
\usepackage{tikz}
\usepackage{tikz-cd}
\usetikzlibrary{arrows,cd,external}
\usepackage{tikz-3dplot}
\usepackage{pgfplots}
\usepackage[exponent-product=\cdot,per-mode=power]{siunitx}

\usepackage{subcaption}

\usepackage{nicefrac}

\definecolor{LightGray}{gray}{0.30}
\definecolor{gray}{rgb}{0.5,0.5,0.7}
\definecolor{dblue}{rgb}{0,0,0.7}
\definecolor{dred}{rgb}{0.72,0,0}
\definecolor{dgreen}{rgb}{0,0.5,0}
\definecolor{dmagenta}{rgb}{0.5,0,0.5}
\definecolor{lblu}{rgb}{0.3,0.3,0.7}

\def\Rd{\color{dred}}

\def\Gd{\color{dgreen}}
\def\Md{\color{dmagenta}}

\def\Bd{\color{dblue}}

\def\div{{ \rm div}}
      \def\grad{{\bf \operatorname{grad}}}
      \def\curl{{\bf curl}}

\def\Hcurlc{{\bf H}({\bf curl};\comdom_\mathrm{c})\,}

\def\Hocurl{{\bf H}_0({\bf curl};\comdom)\,}

\def\Hograd{{\bf H}_0({\bf grad};\comdom)\,}

\def\Hogradc{{\bf H}_0({\bf grad};\comdom_{\mathrm{c}})\,}

\def\Hogradi{{\bf H}_0({\bf grad};\comdom_{\mathrm{i}})\,}

\newcommand{\harm}{\widetilde{\boldsymbol{\mathcal{H}}}}
\newcommand{\harmt}{\widetilde{\boldsymbol{h}}}

\newcommand{\Ngrad}{N_{\phi}}\newcommand{\Nharm}{N_{\harmt}}

\def\SC#1{\ifnum `#1<`1 {\Gd \hat {X}^0_h} \else \ifnum  `#1<`2
  \Bd{\hat X^1_h} \else \ifnum  `#1<`3 \Rd{\hat X^2_h} \else \Md{\hat X^3_h} \fi  \fi  \fi}
\def\XC#1{\ifnum `#1<`1 {\Gd {X}^0_h} \else \ifnum  `#1<`2
  \Bd{ X^1_h} \else \ifnum  `#1<`3 \Rd{ X^2_h} \else \Md{ X^3_h} \fi  \fi  \fi}
\def\XCt#1{\ifnum `#1<`1 {\Gd \tilde {X}^0_h} \else \ifnum  `#1<`2
  \Bd{ \tilde X^1_h} \else \ifnum  `#1<`3 \Rd{ \tilde X^2_h} \else \Md{ \tilde X^3_h} \fi  \fi  \fi}

\theoremstyle{thmstyleone}\newtheorem{theorem}{Theorem}\newtheorem{proposition}[theorem]{Proposition}

\theoremstyle{thmstyletwo}\newtheorem{remark}{Remark}

\theoremstyle{thmstylethree}

\begin{document}

\title{Arbitrary order spline representation of cohomology generators for isogeometric analysis of eddy current problems}

\author[1]{\fnm{Bernard} \sur{Kapidani}
}\email{bernard.kapidani@epfl.ch}

\author[2,3]{\fnm{Melina} \sur{Merkel}
}\email{melina.merkel@tu-darmstadt.de}

\author[2,3]{\fnm{Sebastian} \sur{Schöps}
}\email{sebastian.schoeps@tu-darmstadt.de}

\author[4]{\fnm{Rafael} \sur{V\'azquez}
}\email{rafael.vazquez@usc.es}

\affil*[1]{\orgdiv{Chair of Numerical Modelling and Simulation}, \orgname{École Polytechnique Fédérale de Lausanne}, \orgaddress{\postcode{1015}, \street{Lausanne}, \country{Switzerland}}}
\affil*[2]{\orgdiv{Computational Electromagnetics Group}, \orgname{Technische Universität Darmstadt}, \orgaddress{\postcode{64289}, \street{Darmstadt}, \country{Germany}}}
\affil*[3]{\orgdiv{Centre for Computational Engineering}, \orgname{Technische Universität Darmstadt}, \orgaddress{\postcode{64289}, \street{Darmstadt}, \country{Germany}}}
\affil*[4]{\orgdiv{Departamento de Matem\'atica Aplicada}, \orgname{Universidade de Santiago de Compostela}, \orgaddress{\postcode{15782}, \street{Santiago de Compostela}, \country{Spain}}}

\abstract{The eddy current problem has many relevant practical applications in science, ranging from non-destructive testing to magnetic confinement of plasma in fusion reactors. It arises when electrical conductors are immersed in an external time-varying magnetic field operating at frequencies for which electromagnetic wave propagation effects can be neglected.

 Popular formulations of the eddy current problem either use the magnetic vector potential or the magnetic scalar potential. They have individual advantages and disadvantages. One challenge is related to differential geometry: Scalar potential based formulations run into trouble when conductors are present in non-trivial topology, as approximation spaces must be then augmented with generators of the first cohomology group of the non-conducting domain. 

  For all existing algorithms based on lowest order methods it is assumed that the extension of the graph-based algorithms to high-order approximations requires hierarchical bases for the curl-conforming discrete spaces. However, building on insight on de Rham complexes approximation with splines, we will show in the present submission that the hierarchical basis condition is not necessary. Algorithms based on spanning tree techniques can instead be adapted to work on an underlying hexahedral mesh arising from isomorphisms between spline spaces of differential forms and de Rham complexes on an auxiliary control mesh.
} 


\maketitle

\section{Introduction}

In this paper we discuss the isogeometric discretization of electromagnetic problems that can be simulated by considering the low-frequency `eddy current' approximation to Maxwell's equations. According to \cite[Chapter 8]{Bossavit_1998aa} this approximation neglects the term corresponding to Maxwell's displacement current density. When considering the material properties in Ampere's law one sees that this term is negligible in the conductor inasmuch as the ratio $\varepsilon\omega/\sigma$ can be considered small, where $\varepsilon$ is the permittivity, $\omega$ the frequency and $\sigma$ the conductivity. In the case of induction heating at industrial frequencies, this ratio drops to about $5\cdot10^{-16}$, and Bossavit concludes in \cite{Bossavit_1998aa} that one cannot seriously object to its neglect.

Since the 1980s, many formulations have been proposed to solve the eddy current problem, most of them are based on the combination of potentials, i.e., the magnetic vector potential $\mathbf{A}$ (with $\curl\,\mathbf{A} = \mathbf{B}$, the magnetic flux density), the magnetic scalar potential $\phi$ or $\Omega$, the electric scalar potential $\varphi$ and the electric vector potential $\mathbf{T}$ (with $\curl\,\mathbf{T} = \mathbf{J}$, the current density). Bossavit categorizes in \cite{Bossavit_1985aa} those formulations into two families, magnetic and electric, and identifies the paper \cite{Bossavit_1981aa} and the COMPUMAG conference from 1983 as the corresponding origins. All (vector-)potential formulations suffer from uniqueness issues since gradients lie in the nullspace of the $\curl$-operator and furthermore harmonic fields occur on non-contractible domains, therefore various regularization techniques have been developed \cite{Albanese_1988aa,Dular_1995aa,Bossavit_1999aa,Kettunen_1999aa,Clemens_2002aa}. Later, various experiments were set up with the aim of validating the simulation results computed by the various potential-based finite-element implementations, see e.g. \cite{Bossavit_1988ac,Rodger_1990ab}. The quest of finding optimal formulations for specific problems continues until today, for example various mixed formulations, in particular combinations of $\mathbf{A}$ and $\mathbf{H}$, have been recently proposed in the context of the simulation of superconductors \cite{Grilli_2021aa,Bortot_2020aa,Dular_2021aa}.

Formulations based on the magnetic scalar potential for the non-conductive part of the domain, that is, the $\mathbf{H}-\phi$ and the $\mathbf{T}-\Omega$ formulations, are particularly attractive when solving with the finite element method (FEM). Approximating the magnetic field through a scalar potential drastically reduces the number of degrees of freedom in the non-conductive part of the domain. However, these formulations run into trouble when conductors present non-trivial topology, because the discrete approximation spaces must be augmented with generators of the first cohomology group of the non-conducting domain. Starting from the paper by Bossavit and coauthors in \cite{Kettunen_1998aa}, a lot of research was devoted to save both the use of the computationally cheap magnetic scalar potential and the exact fulfillment of Faraday's law of induction.
The first attempts to give a numerical recipe based on cohomological techniques originally arose in \cite{Kotiuga_1987aa} and was expanded upon in subsequent works~\cite{renSplOmegaFormulation2002,dlotkoAutomaticGenerationCuts2009} where discontinuities of the scalar magnetic potential at mesh aligned cuts, were proposed to handle non simply connected conductors. This approach has nevertheless limitations in the case of particular domains (e.g. knotted conductors) and is usually predicated on the expensive solution of a partial differential equation in three space dimensions.
The work in \cite{Pellikka_2013aa,rodriguezConstructionFiniteElement2013,dlotkoLazyCohomologyGenerators2014,Dlotko_2017aa,Dlotko_2019aa} also built on the same underlying idea of working with the mesh topology information but did not resort to cuts, rather instead providing the additional basis functions as linear combinations of a limited numbers of lowest order edge element basis functions.
For high order FEM commonly hierarchical basis functions are used \cite{Lee_2003aa,Schoberl_2005aa} or constructions based on either moments or weights as degrees of freedom \cite{Los-Santos_2023aa}.

An alternative way to achieve high-order polynomial approximation is through isogeometric analysis (IGA), which was initially introduced within the context of solid mechanics \cite{Hughes_2005aa}. Nowadays, it is commonly interpreted as a special case of FEM that employs spline-based functions for both geometry and solution spaces, thus generalizing the traditional polynomial functions. This approach offers two main advantages. Firstly, the geometry can often be precisely represented due to the use of the same basis functions known as non-uniform rational B-splines (NURBS) employed in computer-aided design. Secondly, the solution becomes globally more regular, as inter-element smoothness can be controlled through the utilization of splines.
In a similar historical development as seen in mixed finite elements (e.g., starting from Nédélec \cite{Nedelec_1980aa}, and leading to works like Arnold et al. \cite{Arnold_2010aa,Arnold_2018aa}), the mechanics community initially explored nodal element spline spaces containing square-integrable gradients. Subsequently, driven by applications in fluid mechanics and electromagnetism, efforts were made to develop high-order div-conforming Raviart-Thomas and curl-conforming Nédélec-element-like multi-patch spline spaces in IGA \cite{Buffa_2010aa,Buffa_2011aa}. More recently, trace-spaces corresponding to these were investigated and integrated into an exact discrete de Rham sequence \cite{Buffa_2019ac}. The curl-conforming spline spaces have been applied to eddy current problems with an $\mathbf{A}$-based formulation in \cite{Friedrich_2020aa} and \cite{Friedrich_2020ab}.

This paper investigates the application of IGA for the solution of the $\mathbf{H}-\phi$ and $\mathbf{T}-\Omega$ formulations. Following the continuous case analogy, the discrete sequence for spline spaces implies that, if the insulator domain is contractible, any square-integrable vector field can be decomposed into the sum of a gradient and a curl, which makes the discretization with a scalar magnetic potential trivial. The extension to non-contractible domains has the same issue explained above for FEM: the space of the scalar magnetic potential must be augmented with generators of the first cohomology groups. We will show that the algorithms of \cite{Dlotko_2017aa} can be easily applied to also work for high order splines. The algorithms are in fact applied to low order FEM spaces defined in an auxiliary hexahedral mesh, called the control mesh. The topological properties are preserved thanks to commutative isomorphisms, introduced in \cite{Buffa_2013aa}, between the FEM spaces of the control mesh and the high order spline spaces. This preservation of the topological properties between the two meshes was already exploited to gauge vector potential formulations in magnetostatics in \cite{Kapidani_2022aa}. 

The paper is structured as follows: \cref{sec:maxwell} introduces Maxwell's equations and the relevant model problems in their strong and weak form. In \cref{sec:iga} isogeometric analysis is summarised, the B-spline spaces forming a de Rham sequence are presented, first in a single-patch domain, and then extended to multi-patch domains. \cref{sec:cohomology} constitutes the main contribution, as it introduces the ingredients needed for the adaptation of the graph based computational homology techniques of \cite{Dlotko_2017aa} to high order spline spaces case. Finally numerical results, are shown in \cref{sec:numerics} to support the theoretical findings. Some conclusions are drawn in \cref{sec:conclusions}.

\section{Maxwell Equations for eddy currents}\label{sec:maxwell}
General electromagnetic fields obey Maxwell's equations \cite{Jackson_1998aa}. However, in the low-frequency regime the magnetoquasistatic field approximation is often sufficient, in which one neglects the displacement current density with respect to the other current densities \cite{Dirks_1996aa,Larsson_2007aa}. This yields the so-called eddy current problem which is given here for the time-harmonic case
\begin{align}\label{eq:maxwell}
	\curl\, \mathbf{E} = - i\omega\mathbf{B},
	&&
	\curl\, \mathbf{H} = \mathbf{J},
	&&
	\div\, \mathbf{B}  = 0
\end{align}
with the two material relations
\begin{align}
	\label{eq:magnetic_material_law}
	\mathbf{B} &= \mu \mathbf{H},
	&
	\quad
	\text{or}
	\quad
	&&
	\mathbf{H} &= \nu \mathbf{B},\\
	\label{eq:ohms_law}
	\mathbf{J}_\textrm{c} &= \sigma \mathbf{E},&
	\quad
	\text{or}
	\quad
	&&
	\mathbf{E} &= \rho \mathbf{J}_\textrm{c},\end{align}
where ${\mathbf{E}}$ is the \emph{electric field strength}, ${\mathbf{B}}$ the \emph{magnetic flux density}, ${\mathbf{H}}$ the \emph{magnetic field 
strength} and $\mathbf{J}=\mathbf{J}_\textrm{c}+\mathbf{J}_\mathrm{s}$ the \emph{total electric current density} consisting of the \emph{conduction} and \emph{source current density} ${\mathbf{J}}_\mathrm{c}$ and ${\mathbf{J}}_\mathrm{s}$, respectively. The magnetic material law \eqref{eq:magnetic_material_law} can be expressed in terms of the \emph{permeability} $\mu$ or the \emph{reluctivity} $\nu$, with $\mu=\nu^{-1}$ and Ohm's law \eqref{eq:ohms_law} either using the \emph{conductivity} $\sigma$ or the \emph{resistivity} $\rho$, with $\sigma=\rho^{-1}$. All the fields are vector-valued, i.e., $\comdom\rightarrow \mathbbm{C}^3$ depending on space, where the computational domain $\comdom\subset\mathbbm{R}^3$ shall be bounded with Lipschitz boundary; the relevance of contractibility will be discussed later. Furthermore, we assume a non-overlapping decomposition of the domain into an insulating and a conducting subdomain, i.e.,
$ 
\comdom
=
\comdom_{\mathrm{i}}
\cup
\comdom_{\mathrm{c}}
$. 
The source current density is impressed in the insulating domain, i.e., $\operatorname{supp}(\mathbf{J}_{\textrm{s}})\subset\comdom_{\mathrm{i}}$, and the conductivity vanishes on domains corresponding to insulators, i.e., $\operatorname{supp}(\sigma)=\comdom_{\mathrm{c}}$. We consider, for simplicity, isotropic materials such that $\sigma:\comdom\rightarrow\mathbbm{R}_{\geq0}$ and $\mu:\comdom\rightarrow\mathbbm{R}_{>0}$ are scalars and uniformly bounded. We will denote by $\Gamma = V_{\mathrm{c}} \cap V_{\mathrm{i}}$ the interface between subdomains.

\bigskip

One rarely solves the system \eqref{eq:maxwell}-\eqref{eq:ohms_law} in its original variables, instead one combines the relevant equations into a \emph{formulation} by exploiting \emph{potentials}. There are two families according to Bossavit~\cite{Bossavit_1985aa}, \emph{electric} formulations, e.g.
$\mathbf{A}-\varphi$ and its variants \cite{Biro_1995aa} as well as \emph{magnetic} formulations based on
$\mathbf{H}-\phi$ 
or
$\mathbf{T}-\Omega$ \cite{Webb_1993aa}. The focus of this paper will be on the magnetic formulations, but we will also use an electric formulation for validation.

\bigskip

The $\mathbf{H}-\phi$ formulation is constructed from
\begin{align*}\mathbf{H} - \mathbf{H}_{\mathrm{s}}
	&=
	\begin{cases}
		\mathbf{H}_\mathrm{c}   & \text{in the conductor }\comdom_{\mathrm{c}},\\
		\grad \,\phi + \harm  & \text{in the insulator }\comdom_{\mathrm{i}},\\
	\end{cases}
\end{align*}
where $\curl \, \mathbf{H}_{\mathrm{s}} = \mathbf{J}_\mathrm{s}$, 
the unknown \emph{magnetic field strength} is restricted to the conductor domain 
$\mathbf{H}_\mathrm{c}:\comdom_\mathrm{c}\rightarrow\mathbbm{C}$, and the \emph{magnetic scalar potential} 
$\phi:\comdom_{\mathrm{i}}\rightarrow\mathbbm{C}$ is restricted to the insulator. 
The harmonic field $\harm: \comdom_{\mathrm{i}}\rightarrow\mathbbm{C}$ is an irrotational field that cannot be represented as a gradient, and takes value zero when the insulator domain $\comdom_{\mathrm{i}}$ is contractible. 

On $\comdom_\textrm{c}$ the choice of the unknowns leads to
\begin{align}
	\curl\bigl(\rho\curl \, \mathbf{H}_{\mathrm{c}}\bigr) + i\omega\mu\mathbf{H}_{\mathrm{c}}&=i\omega\mu\mathbf{H}_{\textrm{s}}, \label{eq:Hphi_in_Vc}\\
\intertext{and on $\comdom_\textrm{i}$ to }
	\div(\mu \mathbf{H}) = \div\bigl(\mu(\grad \phi + \harm)\bigr) +  \div(\mu\mathbf{H}_{\textrm{s}}) = 0, \notag
\end{align}
which are completed with the interface and boundary conditions 
\begin{align} 
\mathbf{H}_c \times \mathbf{n} = (\grad \phi +\harm) \times \mathbf{n} \quad \text{ on } \Gamma, \label{eq:interface_Hphi} \\
	\phi = 0, \qquad \harm \times \mathbf{n} = \mathbf{0} \qquad \text{ on } \partial \comdom, \notag
\end{align}
 where $\mathbf{n}$ is the unit normal vector.
Putting together the previous equations, we aim at weak solutions. By standard procedures the corresponding weak forms are: find $\mathbf{H}_{\mathrm{c}}\in\Hcurlc$,  $\phi\in\Hogradi$ and $\harm \in \widetilde{H}(\comdom_{\mathrm{i}})$, such that 
\begin{align*}
	\int_{\comdom_{\mathrm{c}}} \rho\; \curl \,\mathbf{H}_{\mathrm{c}} \cdot \curl \,\mathbf{w} \,\mathrm{d} \comdom
	+
	\int_{\comdom_{\mathrm{c}}} i\omega\mu \mathbf{H}_{\mathrm{c}}\cdot \mathbf{w}\,\mathrm{d} \comdom
	&=
	-
	\int_{\comdom_{\mathrm{c}}} i\omega\mu \mathbf{H}_{\mathrm{s}}\cdot \mathbf{w}\,\mathrm{d} \comdom \\
\int_{\comdom_{\mathrm{i}}} \mu(\grad \,{\phi} + \harm)\cdot \grad \,w\,\mathrm{d} \comdom
	&= - \int_{\comdom_{\mathrm{i}}} \mu\mathbf{H}_{\mathrm{s}} \cdot \grad \,w\,\mathrm{d} \comdom \\
\int_{\comdom_{\mathrm{i}}} \mu(\grad \,{\phi} + \harm)\cdot \harmt\,\mathrm{d} \comdom
&= - \int_{\comdom_{\mathrm{i}}} \mu\mathbf{H}_{\mathrm{s}} \cdot \harmt\,\mathrm{d} \comdom
\end{align*}
for all $\mathbf{w}\in\Hcurlc$, $w\in\Hogradi$ and $\harmt \in \widetilde{H}(\comdom_{\mathrm{i}})$ that satisfy an interface condition analogous to \eqref{eq:interface_Hphi}, and where the space of harmonic functions $\widetilde{H}(\comdom_{\mathrm{i}})$ is defined as
\begin{equation*}
	\widetilde{H}(\comdom_{\mathrm{i}}) = \left\{ \harmt \in \mathbf{L}^2(\comdom_{\mathrm{i}}) : \curl \, \harmt = \mathbf{0}, \div \, \harmt = 0, \harmt \times \mathbf{n} = \mathbf{0} \text{ on } \partial V\right\}.
\end{equation*} 

We note that the conductivity in the insulator domain is often set to an artificially small positive value, such that one can directly solve \eqref{eq:Hphi_in_Vc} in the whole domain, obtaining a 
pure $\mathbf{H}$-formulation with a single unknown field, at the cost of having more degrees of freedom in the insulator \cite{Dlotko_2019aa}.

\bigskip

The $\mathbf{T}-\Omega$ formulation describes the electromagnetic fields in terms of
\begin{align*}\mathbf{J}_\mathrm{c} &= \curl \,\mathbf{T} &&
	\text{and} &
	\mathbf{H} - \mathbf{H}_{\mathrm{s}} &= 
	\begin{cases}
		\mathbf{T} - \grad \,\Omega, & \text{in the conductor }\comdom_{\mathrm{c}},\\
		- \grad \,\Omega - \harm, & \text{in the insulator }\comdom_{\mathrm{i}},
	\end{cases}
\end{align*}
with 
the \emph{electric vector potential}
$\mathbf{T}:\comdom_{\mathrm{c}}\rightarrow\mathbbm{C}^{3}$,
the (globally defined) \emph{magnetic scalar potential} 
$\Omega:\comdom\rightarrow\mathbbm{C}$, and the harmonic field $\harm: \comdom_{\mathrm{i}}\rightarrow\mathbbm{C}$ defined as above, and which is nonzero only in the case of a non-contractible domain $\comdom_{\mathrm{i}}$. 
In the conductor $\comdom_\textrm{c}$ the system reads
\begin{align*}
	\curl\bigl(\rho\curl \,\mathbf{T} \bigr) + i\omega\bigl(\mu\mathbf{T}-\mu\grad \,\Omega \bigr)&=-i\omega\mu\mathbf{H}_{\textrm{s}}\\
	\div\bigl(\mu\mathbf{T}-\mu\grad \,\Omega\bigr) &= -\div(\mu\mathbf{H}_{\textrm{s}}),
\intertext{and on $\comdom_\textrm{i}$}
	\div (\mu \mathbf{H}) = \div(\mu\mathbf{H}_{\textrm{s}}) - \div\bigl(\mu(\grad \, \Omega + \harm)\bigr) = 0.
\end{align*}
The interface and boundary conditions of the problem are
\begin{align}\label{eq:interface_Tomega}
	\mathbf{T} \times \mathbf{n} = - \harm \times \mathbf{n} \quad \text{ on } \Gamma, \\
	\Omega = 0, \qquad \harm \times \mathbf{n} =  \mathbf{0}, \text{ on } \partial V.
\end{align}
The problem in weak formulation is then given by: find $\mathbf{T}\in\Hcurlc$, $\Omega\in\Hograd$ and $\harm \in \widetilde{H}(\comdom_{\mathrm{i}})$, satisfying the interface condition \cref{eq:interface_Tomega}, such that
\begin{align*}
	\int_{\comdom_{\mathrm{c}}} \rho\; \curl \, \mathbf{T} \cdot \curl \,\mathbf{w} \,\mathrm{d} \comdom
	+
	\int_{\comdom_{\mathrm{c}}} i\omega\mu \mathbf{T}\cdot \mathbf{w}\,\mathrm{d} \comdom
	-
	\int_{\comdom_{\mathrm{c}}} i\omega\mu \grad \,\Omega\cdot \mathbf{w}\,\mathrm{d} \comdom
	&=
	-
	\int_{\comdom_{\mathrm{c}}} i\omega\mu \mathbf{H}_{\mathrm{s}}\cdot \mathbf{w}\,\mathrm{d} \comdom
	\\
- \int_{\comdom_{\mathrm{c}}} \mu\mathbf{T} \cdot \grad \,w \,\mathrm{d} \comdom
	+
	\int_{\comdom} \mu\grad \,\Omega \cdot \grad \,w \,\mathrm{d} \comdom
	+ 
	\int_{\comdom_\mathrm{i}} \mu \harm \cdot \grad \,w \,\mathrm{d} \comdom
	&= 
	\int_{\comdom} \mu\mathbf{H}_{\mathrm{s}} \cdot \grad \,w\,\mathrm{d} \comdom \\
\int_{\comdom_{\mathrm{i}}} \mu(\grad \,\Omega + \harm)\cdot \harmt\,\mathrm{d} \comdom
	&= 
	\int_{\comdom_{\mathrm{i}}} \mu\mathbf{H}_{\mathrm{s}} \cdot \harmt\,\mathrm{d} \comdom
\end{align*}
for all $\mathbf{w}\in\Hcurlc$, $w\in\Hograd$ and $\harmt \in \widetilde{H}(\comdom_{\mathrm{i}})$ that satisfy an interface condition analogous to \eqref{eq:interface_Tomega}.
 To ensure uniqueness of the solution in the conductor, an additional gauging condition for $\mathbf{T}$ must be considered \cite{Webb_1993aa}.

Note that the $\mathbf{H}-\phi$ and $\mathbf{T}-\Omega$ formulations are very similar, the difference is only in the fact where the unknowns are supported. The resulting advantage of the (more complex) $\mathbf{T}-\Omega$ formulation is that it remains stable in the frequency limit $\omega\to0$ \cite{Webb_1993aa}. However, both formulations have issues for non-contractible domains, in particular for the approximation of the harmonic field $\harm$ in a complex domain, as will be discussed in \cref{sec:cohomology}.

Finally, the $\mathbf{A}-\varphi$ formulation \cite{Kameari_1990aa,Biro_1989aa,Bossavit_1998aa,Manges_1997aa,Clemens_2002aa} that we will use for results comparison is based on
\begin{align*}\mathbf{B}&=\curl \,\mathbf{A}
	&& \text{and} &
	\mathbf{E}&=-i\omega\mathbf{A} -\grad \,\varphi
\end{align*}
involving the \emph{magnetic vector potential} 
$\mathbf{A}:\comdom\rightarrow\mathbbm{C}^{3}$
and the \emph{electric scalar potential} 
$\varphi:\comdom_\textrm{c}\rightarrow\mathbbm{C}$. The system of equations is then given on 
$\comdom_\textrm{c}$ by \begin{align}
	\label{eq:a-phi-problem1}
	\curl\bigl(\nu\curl \,\mathbf{A}\bigr) + i\omega\sigma\mathbf{A} + \sigma\grad \, \varphi&=0
\\
	\label{eq:a-phi-problem2}
	\div\bigl(i\omega\sigma\mathbf{A}-\sigma\grad \,\varphi \bigr) &= 0
	\intertext{and on $\comdom_\textrm{i}$ by}
	\label{eq:a-phi-problem3}
	\curl\bigl(\nu\curl \,\mathbf{A}\bigr)  &=\mathbf{J}_{\mathrm{s}}.
\end{align}
We apply the following Dirichlet interface and boundary conditions
\begin{align*}
	\mathbf{A} \times \mathbf{n} = \mathbf{0} \qquad & \text{ on } \partial \comdom, \\
	\varphi = 0 \qquad & \text{ on } \Gamma,
\end{align*}
see \cite[System (45)-(55)]{Biro_1989aa}. Let us formulate the weak form for \eqref{eq:a-phi-problem1}-\eqref{eq:a-phi-problem3} directly on both domains: find $\mathbf{A}\in\Hocurl$ and $\varphi\in\Hogradc$ such that
\begin{align*}
	\!\int_{\comdom}\! \nu\curl \,\mathbf{A} \cdot \curl \,\mathbf{w} \,\mathrm{d}\comdom
	+
\!\int_{\comdom_{\mathrm{c}}}\!i\omega\sigma\mathbf{A}\cdot{\mathbf{w}}\,\mathrm{d} \comdom
	+
	\!\int_{\comdom_{\mathrm{c}}}\!\sigma\mathbf{w}\cdot\grad{\varphi}\,\mathrm{d} \comdom
	&=
	\!\int_{\comdom_{\mathrm{i}}}\!\mathbf{J}_{\mathrm{s}}\cdot{\mathbf{w}}\,\mathrm{d} \comdom
	\\
	\!\int_{\comdom_{\mathrm{c}}}\!i \omega \sigma\grad\, w \cdot \mathbf{A} \,\mathrm{d}\comdom
	+
	\!\int_{\comdom_{\mathrm{c}}}\! \sigma\grad \,w \cdot \grad \,\varphi \,\mathrm{d}\comdom
	&=
	0
\end{align*}
for all $\mathbf{w}\in\Hocurl$ and $w\in\Hogradc$. In general, the magnetic flux density $\mathbf{B}$ defines the magnetic vector potential $\mathbf{A}$ only up to an irrotational field. For a unique solution an additional gauging condition is required, and one could use Coulomb gauging or tree-cotree gauging, for instance. Another possibility is to add an artificially small conductivity in the region $\comdom_{\mathrm{i}}$, that also guarantees uniqueness of the solution.

\bigskip

\section{Isogeometric Analysis}\label{sec:iga}
In the first part of this section we recall the construction of tensor-product spline spaces that satisfy a de Rham diagram, and their relation through commuting isomorphisms with lowest order finite elements on a structured hexahedral grid \cite{Buffa_2011aa,Buffa_2013aa}. In the second part, we show how the definitions are extended to the multi-patch case. Finally, we apply the discrete spaces of the sequence to the discretization of the variational formulations of Section~\ref{sec:maxwell}.

\subsection{Spline complex in single-patch domains}
B-spline basis functions of degree $p$ are defined from an ordered knot vector $\Xi = \{\xi_1, \ldots, \xi_{n+p+1}\}$, with $0 \le \xi_i \le \xi_{i+1} \le 1$ for every $i$. We assume that the knot vector is open, which means that the first and last knot are repeated exactly $p+1$ times. We denote by $\hat B_i^p$, for $i=1, \ldots, n$, the $i$-th basis function of degree $p$ in the reference domain $(0,1)$, that can be computed with the Cox-de Boor formula \cite[Chapter IX]{de-Boor_2001aa}, and by $S_p(\Xi)$ the space they span. From the knot vector $\Xi$, we also define the modified knot vector $\Xi'$ by removing the first and last knot, and the spline space $S_{p-1}(\Xi')$ of degree $p-1$. For this space, we replace the standard B-spline basis by the Curry-Schoenberg B-splines $\hat D_i^{p-1} = \frac{p}{\xi_{i+p}-\xi_i} \hat{B}_{i}^{p-1}$, as it was first done in \cite{Ratnani_2012aa}, from which we get the following expression for the derivative of B-splines: 
\begin{equation} \label{eq:derivative}
\frac{d \hat{B}_i^p}{d x} = \hat{D}_{i-1}^{p-1} - \hat{D}_{i}^{p-1}, \quad \text{ for } i = 1, \ldots, n,
\end{equation}
where we assume the convention $D_0^{p-1} = D_{n}^{p-1} \equiv 0$. This particular choice of the basis functions allows to write the derivative as an incidence matrix $\mathbb{G}_n \in \mathbb{R}^{(n-1) \times n}$, which takes the form
\begin{equation} \label{eq:incidence}
\mathbb{G}_n = 
\begin{bmatrix}
-1 & 1 \\
& \ddots & \ddots \\
& & -1 & 1
\end{bmatrix}.
\end{equation}

In the reference domain $\hat{V} = (0,1)^3$, multivariate B-splines are defined by tensor-product of univariate B-splines. Given a vector degree $\mathbf{p} = (p_1, p_2, p_3)$, and the knot vectors $\Xi_j$, for $j=1,2,3$, the multivariate basis functions are given by
\[
\hat B_{\mathbf{i}}^{\mathbf{p}} (\boldsymbol{\xi}) = \hat B_{i_1}^{p_1}(\xi_1)  \, \hat B_{i_2}^{p_2}(\xi_2) \, \hat B_{i_3}^{p_3}(\xi_3), \quad \text{ for } \boldsymbol{\xi} = (\xi_1, \xi_2, \xi_3) \in \hat {V},
\]
where the multi-index $\mathbf{i} = (i_1, i_2, i_3)$ satisfies $1 \le i_k \le n_k$, with $n_k$ the number of basis functions in each direction. The multivariate B-splines span the spline space $S_{p_1,p_2,p_3}(\Xi_1,\Xi_2,\Xi_3) = S_{p_1}(\Xi_1) \otimes S_{p_2}(\Xi_2) \otimes S_{p_3}(\Xi_3)$. For simplicity we assume from now on that the degree is equal in all directions, that is, $p_j = p$ for every $j$. A discrete de Rham diagram of spline spaces can be defined as in \cite{Buffa_2010aa,Buffa_2011aa}, and takes the form
\begin{equation*}
\begin{tikzcd}
S_p^0(\hat V) \arrow[r,"\grad"] & S_p^1(\hat V) \arrow[r,"\curl"] & S_p^2(\hat V) \arrow[r,"\div"] & S_p^3(\hat V),
\end{tikzcd}\label{eq:spline_derham}
\end{equation*}
where the spaces are obtained by combining splines of degrees $p$ and $p-1$, in particular
\begin{align*}
& S_p^0(\hat V) = S_{p,p,p}(\Xi_1,\Xi_2,\Xi_3), \\
& S_p^1(\hat V) = S_{p-1,p,p}(\Xi'_1,\Xi_2,\Xi_3) \times S_{p,p-1,p}(\Xi_1,\Xi'_2,\Xi_3) \times S_{p,p,p-1}(\Xi_1,\Xi_2,\Xi'_3), \\
& S_p^2(\hat V) = S_{p,p-1,p-1}(\Xi_1,\Xi'_2,\Xi'_3) \times S_{p-1,p,p-1}(\Xi'_1,\Xi_2,\Xi'_3) \times S_{p-1,p-1,p}(\Xi'_1,\Xi'_2,\Xi_3), \\
& S_p^3(\hat V) = S_{p-1,p-1,p-1}(\Xi'_1,\Xi'_2,\Xi'_3).
\end{align*}

A spline geometry is built as a map from the reference domain $\hat{V}$ to the physical domain $V \subset \mathbb{R}^3$ by associating a control point $\mathbf{P}_{\mathbf{i}} \in \mathbb{R}^3$ to each basis function, namely
\begin{equation} \label{eq:param}
\mathbf{F}(\boldsymbol{\xi}) = \sum_{\mathbf{i}} {\mathbf{P}}_{\mathbf{i}} \hat B_{\mathbf{i}}^{\mathbf{p}} (\boldsymbol{\xi}), \quad \boldsymbol{\xi} \in \hat{V}.
\end{equation}
It is also possible to use a NURBS geometry, instead of a spline geometry, replacing in the previous expression the B-splines by NURBS basis functions, which are defined by associating a weight to each B-spline, see \cite[Chapter~2]{Hughes_2005aa} for details. We can then construct a discrete de Rham sequence in the physical domain $V$ as
\begin{equation}
\begin{tikzcd}
S_p^0(V) \arrow[r,"\grad"] & S_p^1(V) \arrow[r,"\curl"] & S_p^2(V) \arrow[r,"\div"] & S_p^3(V),
\end{tikzcd}\label{eq:spline_derham_phys}
\end{equation}
where each space is obtained by applying a suitable pull-back to the corresponding space in the reference domain, and which is based on the map $\mathbf{F}$, see \cite{Buffa_2011aa}.

\begin{figure}[t]
\centering
\includegraphics[trim=2cm 2cm 1cm 2cm, clip, width=0.4\textwidth]{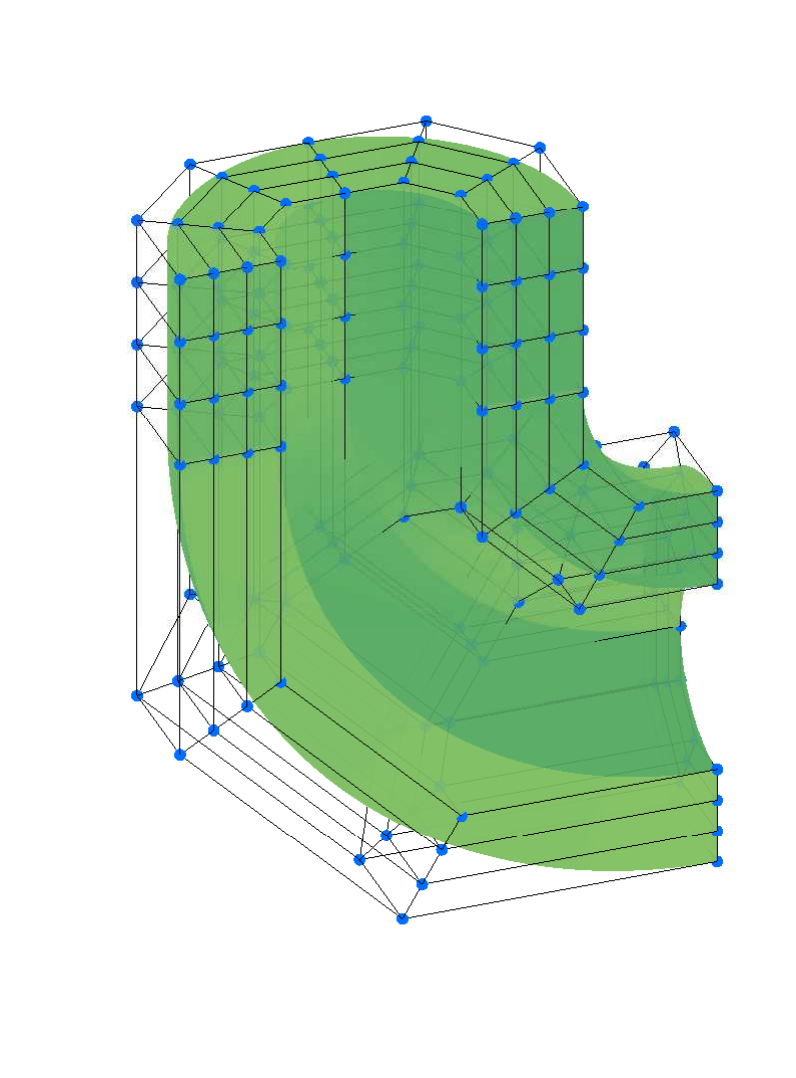}
\caption{An example of a parametrized domain $V$ (in green), and the corresponding control points (blue dots), which determine the control mesh.}\label{fig:control_mesh}
\end{figure}
The set of control points in \eqref{eq:param} determines a structured grid of hexahedral trilinear elements, called the control mesh, which encloses a polyhedral domain that we denote by $V_C$. An example of a parametrized geometry, and the corresponding control mesh, is given in Figure~\ref{fig:control_mesh}. It is clear, from the definition of $\mathbf{F}$, that each basis function of $S_p^0(V)$ is associated to a control point, that is, to a vertex of the control mesh. In an analogous way, the basis functions of the spaces $S_p^1(V), S_p^2(V)$ and $S_p^3(V)$ are respectively associated to the edges, faces and cells of the control mesh. Let us now define, on the hexahedral control mesh, the de Rham sequence of lowest order finite element spaces, that we will denote by $Z_h^k(V_C)$ for $k=0, \ldots, 3$. In particular, $Z_h^0(V_C)$ is the space spanned by the standard nodal finite elements, $Z_h^1(V_C)$ is the space of edge elements, $Z_h^2(V_C)$ is the space of lowest order face elements, and $Z_h^3(V_C)$ is the space of piecewise constants, all of them defined on the control mesh. Since the basis functions of these spaces are associated with geometrical entities of the control mesh, and we have seen that the same is true for the spline basis functions, we can automatically define isomorphisms $I_h^k : S_p^k(V) \rightarrow Z_h^k(V_C)$ that take each spline basis function to the finite element basis function associated to the same geometrical entity. Thanks to the choice of the basis functions and the expression of the derivative in \eqref{eq:derivative}, these isomorphisms commute with the differential operators, giving the following commutative diagram \cite{Buffa_2013aa}:
\begin{equation} \label{eq:isomorphisms}
\begin{tikzcd}
S_p^0(V) \arrow[r,"\grad"] \arrow[d,"I_h^0"] & S_p^1(V) \arrow[r,"\curl"] \arrow[d,"I_h^1"] & S_p^2(V) \arrow[r,"\div"] \arrow[d,"I_h^2"] & S_p^3(V) \arrow[d,"I_h^3"] \\
Z_h^0(V_C) \arrow[r,"\grad"] & Z_h^1(V_C) \arrow[r,"\curl"] & Z_h^2(V_C) \arrow[r,"\div"] & Z_h^3(V_C).
\end{tikzcd}
\end{equation}
These commutative isomorphisms will play a key role in the construction of generators of the cohomology groups in Section~\ref{sec:cohomology}.

Finally, another interesting property which stems from the choice of the basis functions is that the gradient, curl and divergence operators can be expressed as incidence matrices of the control mesh, equivalently to what is done for instance in the finite integration technique (FIT) \cite{Clemens_2001aa} or lowest order Whitney elements \cite[Section 5.2]{Bossavit_1998aa}. More precisely, the gradient operator is given by the vertex-edge incidence matrix $\mathbb{G}$, the curl is given by the edge-face incidence matrix $\mathbb{C}$, and the divergence is given by the face-cell incidence matrix $\mathbb{D}$ for the control mesh. 
Thanks to the tensor-product structure of the spaces, and the fact that the control mesh is a structured mesh, these matrices can be computed from Kronecker tensor-products as \cite{Holderied_2021aa}
\[
\mathbb{G} =
\begin{bmatrix}
\mathbb{I}_{n_3} \otimes \mathbb{I}_{n_2} \otimes \mathbb{G}_{n_1} \\
\mathbb{I}_{n_3} \otimes \mathbb{G}_{n_2} \otimes \mathbb{I}_{n_1} \\
\mathbb{G}_{n_3} \otimes \mathbb{I}_{n_2} \otimes \mathbb{I}_{n_1}
\end{bmatrix}, \qquad
\mathbb{D} =
\begin{bmatrix}
\mathbb{I}_{n_3 -1} \otimes \mathbb{I}_{n_2 -1} \otimes \mathbb{G}_{n_1} \\
\mathbb{I}_{n_3 -1} \otimes \mathbb{G}_{n_2} \otimes \mathbb{I}_{n_1 -1} \\
\mathbb{G}_{n_3} \otimes \mathbb{I}_{n_2 -1} \otimes \mathbb{I}_{n_1 -1}
\end{bmatrix}^\top,
\]
\[
\mathbb{C} =
\begin{bmatrix}
0 &
- \mathbb{G}_{n_3} \otimes \mathbb{I}_{n_2 -1} \otimes \mathbb{I}_{n_1} &
\mathbb{I}_{n_3-1} \otimes \mathbb{G}_{n_2} \otimes \mathbb{I}_{n_1} \\
\mathbb{G}_{n_3} \otimes \mathbb{I}_{n_2} \otimes \mathbb{I}_{n_1 -1} &
0 &
- \mathbb{I}_{n_3 -1} \otimes \mathbb{I}_{n_2} \otimes \mathbb{G}_{n_1} \\
- \mathbb{I}_{n_3} \otimes \mathbb{G}_{n_2} \otimes \mathbb{I}_{n_1 -1} &
\mathbb{I}_{n_3} \otimes \mathbb{I}_{n_2 -1} \otimes \mathbb{G}_{n_1} &
0 
\end{bmatrix},
\]
where $\mathbb{I}_n$ is an identity matrix of size $n \times n$, and the $\mathbb{G}_{n_k}$ matrices correspond to the matrices as in \eqref{eq:incidence} for the knot vectors in the three Cartesian directions.
It is also worth noting that, as in FIT or lowest order finite elements, the incidence matrices can be exploited for the assembly of the matrices arising from the discretization of the weak formulations in Section~\ref{sec:maxwell}. For instance, denoting by $\mathbb{M}^k$ the mass matrix associated to $S_p^k(V)$, the $k$th space in the de Rham diagram, the grad-grad and the curl-curl operators can be respectively computed using the matrices
\[
\mathbb{A}^0 = \mathbb{G}^\top \mathbb{M}^1 \mathbb{G}, \qquad
\mathbb{A}^1 = \mathbb{C}^\top \mathbb{M}^2 \mathbb{C}.
\]

\begin{remark}
We have assumed for simplicity that the map $\mathbf{F}$ uses the same functions as the first space of the de Rham diagram. However, in practice the kind of mappings that can be used to construct the geometry is much wider, and not necessarily restricted to spline or NURBS spaces, see \cite{Buffa_2011aa} or \cite[Section~3]{Beirao-da-Veiga_2014aa} for the theoretical assumptions on $\mathbf{F}$. In this case, the existence of the control mesh is not necessarily guaranteed. A hexahedral mesh with the same structure can be obtained replacing the control points by the image through $\mathbf{F}$ of the Greville points (or knot averages), and the incidence matrices, the lowest order finite element spaces and the isomorphisms can be defined using this Greville mesh \cite{Buffa_2013aa}. For simplicity, from now on we will keep referring to the control mesh.
\end{remark}

\begin{remark}
	To impose homogeneous Dirichlet boundary conditions in the IGA setting, one simply has to remove the degrees of freedom of control variables associated to the boundary. For instance, for the space $S_p^0(\comdom) \cap \Hograd$ the boundary control points are not considered, while for $S_p^1(\comdom) \cap \Hocurl$ the boundary edges are not considered. To simplify notation, in the following we will use the same notation for discrete spaces with or without boundary conditions.
\end{remark}

\subsection{Spline complex in multi-patch domains}
Although the use of a mapping as in \eqref{eq:param} allows to work on curved domains, it is restricted to simple geometries which are the image of a single cube, what is usually called a patch. More complex domains can be defined as the images of several patches, in the form $\overline V = \bigcup_{j=1}^N \overline{V^{(j)}}$, and assuming that each subdomain is obtained from a mapping $\mathbf{F}^{(j)} : (0,1)^3 \rightarrow V^{(j)}$ as in \eqref{eq:param}. To construct a basis of the multi-patch domain we assume that the subdomain partition is conforming, in the sense that the mappings and the knot vectors coincide, up to an affine transformation, on every interface between patches, see \cite{Kleiss_2012ab} and \cite[Assumption~3.5]{Beirao-da-Veiga_2014aa} for a rigorous description of the assumptions. This implies that on the interface the meshes also coincide, and there is a one-to-one correspondence of the basis functions from the two involved patches.

\begin{figure}[t]
\centering
\includegraphics[trim=1cm 2cm 1cm 1cm, clip, width=0.42\textwidth]{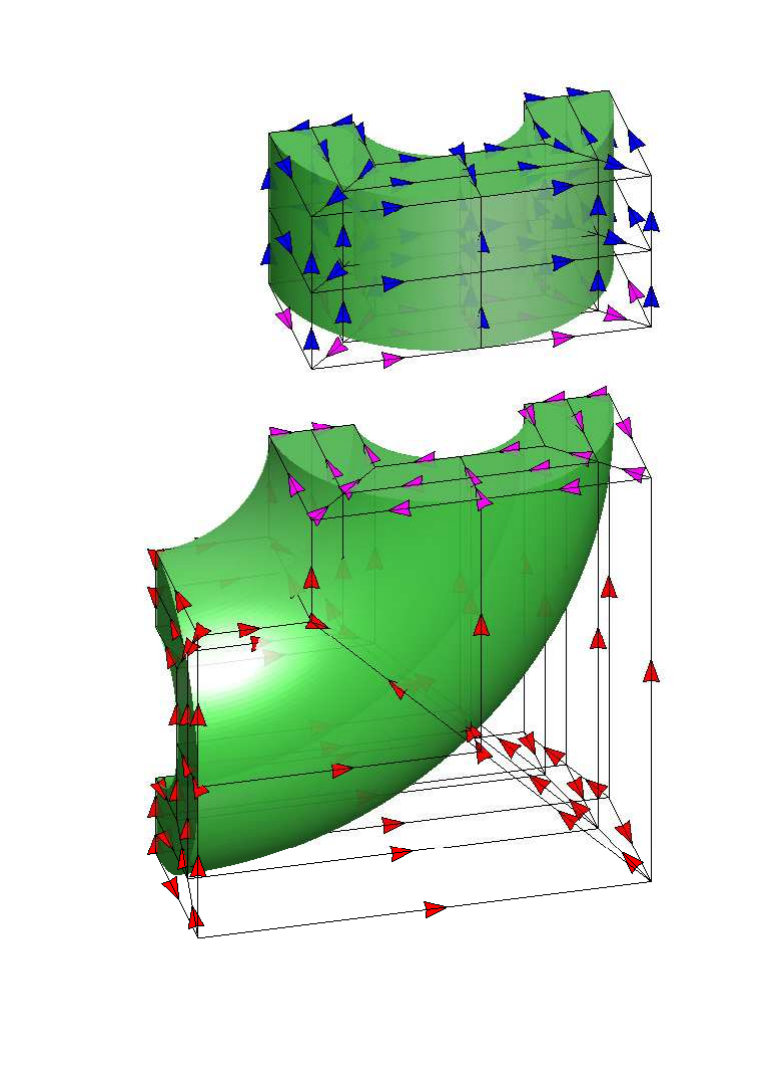}
\includegraphics[trim=1cm 2cm 1cm 1cm, clip, width=0.35\textwidth]{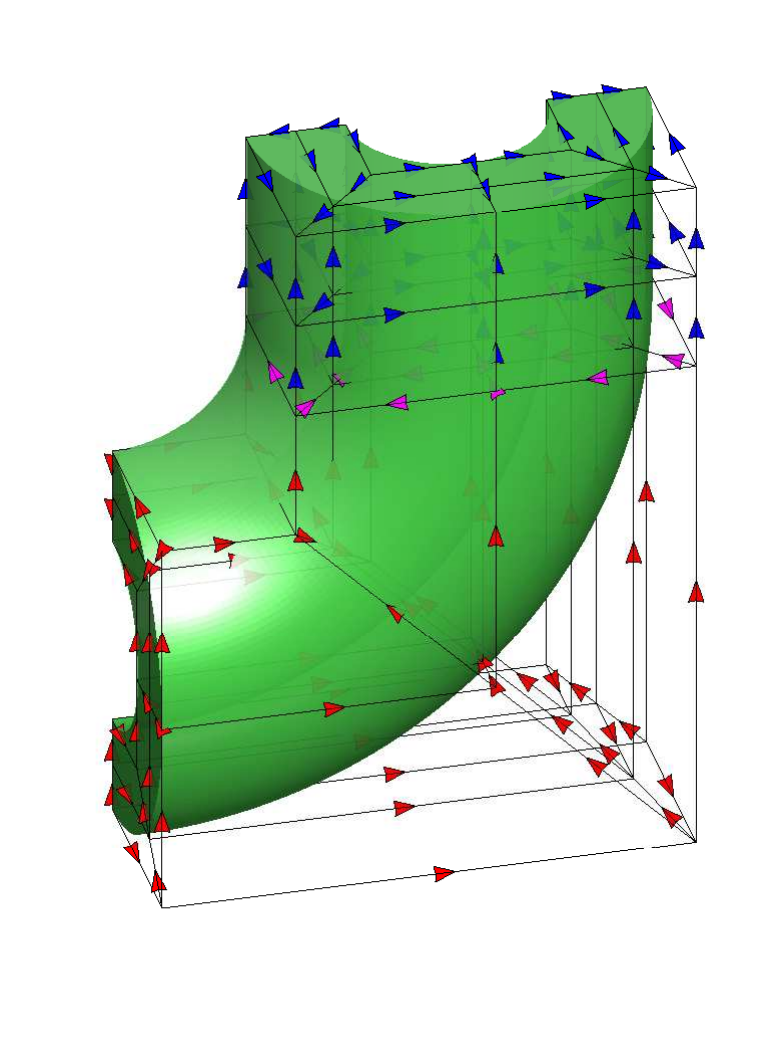}
\caption{The process of gluing together the basis functions of $S_p^1(V)$ for a two-patch domain. On the left, the basis functions of the two patches are represented with their local orientation, induced by the one in the parametric domain. On the right, the magenta functions on the interface have been glued together, and one of the two orientations had to be chosen.} \label{fig:multipatch}
\end{figure}

To define the spline spaces in the multi-patch domain, we first introduce a sequence of spline spaces on each patch $V^{(j)}$, denoted by $S_p^k(V^{(j)})$ for $k=0, \ldots, 3$, as in \eqref{eq:spline_derham_phys}. Then, with harmless abuse of notation, the spaces on the multi-patch domain are defined as
\begin{align*}
& S_p^0(V) = \{ u \in H^1(V) : u|_{V^{(j)}} \in S_p^0(V^{(j)}) \}, \\
& S_p^1(V) = \{ \mathbf{u} \in \mathbf{H}(\mathbf{curl}, V): \mathbf{u}|_{V^{(j)}} \in S_p^1(V^{(j)}) \}, \\
& S_p^2(V) = \{ \mathbf{u} \in \mathbf{H}(\mathrm{div}, V): \mathbf{u}|_{V^{(j)}} \in S_p^2(V^{(j)}) \}, \\
& S_p^3(V) = \{ u \in L^2(V) : u|_{V^{(j)}} \in S_p^3(V^{(j)}) \}.
\end{align*}
In practice, the conformity condition (such as $\mathbf{u} \in \mathbf{H}(\mathbf{curl}; V)$ for $S_p^1(V)$) is imposed by gluing together the basis functions on the interfaces, thanks to the one-to-one correspondence implied by the assumptions above. This can be easily done with the help of the control mesh. In fact, we can introduce a control mesh for each patch exactly as in the single-patch case, and thanks to the assumptions on the subdomain partition, the control points from each patch on the interface must coincide. Therefore, we can glue functions together by identifying the control points, control edges and control faces on the interfaces, in an analogous way as in finite elements local degrees of freedom are glued between elements leading in both cases to global $C^0$-continuity. In this process, the orientation of the vector-valued basis functions of $S_p^1(V)$ and $S_p^2(V)$ must be taken into account, analogously to the definition of the orientation for edge and face finite elements. See Figure~\ref{fig:multipatch} for an example of how to glue together the functions of two patches for the space $S^1_p(V)$.

We have just mentioned that the control meshes from different patches have coinciding control points at the interfaces. Therefore, we can define a global control mesh, which is locally structured on each patch but globally unstructured. As in the single-patch case, but now up to nontrivial cohomology, we can construct the spaces of lowest order finite elements on this control mesh, which we still denote by $Z_h^k(V_C)$, and define commutative isomorphisms between the multi-patch spline spaces and the lowest order finite element spaces on the control mesh, in such a way that the commutative diagram \eqref{eq:isomorphisms} remains valid in the multi-patch setting.

Furthermore, in similar fashion to the three-dimensional spaces, once the interface manifold $\Gamma$ between the $V_{\mathrm{i}}$ and $V_{\mathrm{c}}$ subdomains has been discretized by a multipatch geometry $\bigcup_{0<j<N} \Gamma^{(j)}$, the results of \cite{Buffa_2019ac} can be used to define the trace spaces $S_p^0(\Gamma)$, $S_p^1(\Gamma)$, $S_p^2(\Gamma)$ 
which also constitute an exact sequence up to non trivial cohomology:
\begin{equation}
	\begin{tikzcd}
	S_p^0(\Gamma) \arrow[r,"\grad_\Gamma"] & S_p^1(\Gamma) \arrow[r,"\mathrm{curl}_\Gamma"] & S_p^2(\Gamma),
	\end{tikzcd}
	\label{eq:spline_derham_gamma}
\end{equation}
where we note that surface gradients and curl operators are used.

It is finally worth noting that the expression of the differential operators in terms of the incidence matrices of the control mesh also remains valid. Since the global mesh is now unstructured, the Kronecker product expressions of the incidence matrices are not valid anymore, but the tensor-product structure can be exploited by computing the incidence matrices patchwise. Similarly, the assembly of the matrices in the multi-patch setting can also be done with the help of the incidence matrices.

\subsection{Discrete weak formulation} \label{subsec:weak_IGA}
The solution with IGA of the $\mathbf{H}-\phi$ or the $\mathbf{T}-\Omega$ formulations in the case of contractible domains does not pose any particular problem: the vector fields $\mathbf{H}$ or $\mathbf{T}$ are approximated with curl-conforming splines, while the scalar magnetic potential $\phi$ or $\Omega$ is approximated with gradient-conforming splines, that is, the first space in the de Rham sequence.

The main challenge for both formulations is to find a discrete version of the field $\harm$ when the domain is non-contractible. Obviously, all gradient fields are irrotational fields, also in the discrete case. But thanks to the existence of commutative projectors between the continuous and the discrete de Rham sequences \cite{Buffa_2011aa}, there are other irrotational functions which are not given by gradients. This can be written in the form
\begin{equation*}
\{\mathbf{u}_h \in S_p^1(\comdom_{\mathrm{i}}) : \curl \, \mathbf{u_h} = \mathbf{0} \} = \{\grad \, \psi_h : \psi_h \in S_p^0(\comdom_{\mathrm{i}})\} + W_h(V_{\mathrm{i}}),
\end{equation*}
where $W_h(\comdom_{\mathrm{i}}) \subset S_p^1(\comdom_{\mathrm{i}})$ satisfies $\dim{W_h(\comdom_{\mathrm{i}})} = \dim{\widetilde{H}(\comdom_\mathrm{i})}$. While the functions in $W_h(\comdom_{\mathrm{i}})$ are not necessarily harmonic (i.e., it may happen that $W_h(\comdom_{\mathrm{i}}) \not \subset \widetilde{H}(\comdom_{\mathrm{i}})$), they must be generators of the first cohomology group to represent all discrete irrotational functions. 

Assuming that a discrete space $W_h(\comdom_{\mathrm{i}})$ is known, we can write the discrete version of the $\mathbf{H}-\phi$ formulation as: find $\mathbf{H}_h\in S_p^1(\comdom_{\mathrm{c}})$,  $\phi_h \in S_p^0(\comdom_\mathrm{i})$ and $\harm_h \in W_h(\comdom_{\mathrm{i}}) $, such that they satisfy the interface condition 
\begin{equation} \label{eq:interface_Hphi_disc}
	\mathbf{H}_h \times \mathbf{n} = (\grad \phi_h +\harm_h) \times \mathbf{n} \quad \text{ on } \Gamma,
\end{equation}
and moreover they satisfy
\begin{align*}
	\int_{\comdom_{\mathrm{c}}} \rho\; \curl \,\mathbf{H}_{h} \cdot \curl \,\mathbf{w}_h \,\mathrm{d} \comdom
	+
	\int_{\comdom_{\mathrm{c}}} i\omega\mu \mathbf{H}_{h}\cdot \mathbf{w}_h\,\mathrm{d} \comdom
	&=
	-
	\int_{\comdom_{\mathrm{c}}} i\omega\mu \mathbf{H}_{\mathrm{s}}\cdot \mathbf{w}_h\,\mathrm{d} \comdom \\
\int_{\comdom_{\mathrm{i}}} \mu(\grad \,{\phi}_h + \harm_h)\cdot \grad \,w_h\,\mathrm{d} \comdom
	&= - \int_{\comdom_{\mathrm{i}}} \mu\mathbf{H}_{\mathrm{s}} \cdot \grad \,w_h\,\mathrm{d} \comdom \\
\int_{\comdom_{\mathrm{i}}} \mu(\grad \,{\phi}_h + \harm_h)\cdot \harmt_h\,\mathrm{d} \comdom
	&= - \int_{\comdom_{\mathrm{i}}} \mu\mathbf{H}_{\mathrm{s}} \cdot \harmt_h\,\mathrm{d} \comdom
\end{align*}
for all $\mathbf{w}_h\in S_p^1(\comdom_{\mathrm{c}})$, $w_h\in S_p^0(\comdom_{\mathrm{i}})$ and $\harmt_h \in W_h(\comdom_{\mathrm{i}}) $ that also satisfy the interface condition \eqref{eq:interface_Hphi_disc}.

In a completely analogous way, the discrete version of the $\mathbf{T}-\Omega$ formulation is given by: find $\mathbf{T}_h \in S_p^1(\comdom_{\mathrm{c}})$, $\Omega_h\in S_p^0(\comdom)$ and $\harm_h \in W_h(\comdom_{\mathrm{i}})$, satisfying the interface condition $\mathbf{T}_h \times \mathbf{n} = -\harm_h \times \mathbf{n}$ on the interface $\Gamma$, and such that
\begin{align*}
	\int_{\comdom_{\mathrm{c}}} \rho\; \curl \, \mathbf{T}_h \cdot \curl \,\mathbf{w}_h \,\mathrm{d} \comdom
	+
	\int_{\comdom_{\mathrm{c}}} i\omega\mu \mathbf{T}_h\cdot \mathbf{w}_h\,\mathrm{d} \comdom
	-
	\int_{\comdom_{\mathrm{c}}} i\omega\mu \grad \,\Omega_h\cdot \mathbf{w}_h\,\mathrm{d} \comdom
	&=
	-
	\int_{\comdom_{\mathrm{c}}} i\omega\mu \mathbf{H}_{\mathrm{s}}\cdot \mathbf{w}_h\,\mathrm{d} \comdom
	\\
- \int_{\comdom_{\mathrm{c}}} \mu\mathbf{T}_h \cdot \grad \,w_h \,\mathrm{d} \comdom
	+
	\int_{\comdom} \mu\grad \,\Omega_h \cdot \grad \,w_h \,\mathrm{d} \comdom
	+ 
	\int_{\comdom_\mathrm{i}} \mu \harm_h \cdot \grad \,w_h \,\mathrm{d} \comdom
	&= 
	\int_{\comdom} \mu\mathbf{H}_{\mathrm{s}} \cdot \grad \,w_h\,\mathrm{d} \comdom \\
\int_{\comdom_{\mathrm{i}}} \mu(\grad \,\Omega_h + \harm_h)\cdot \harmt_h\,\mathrm{d} \comdom
	&= 
	\int_{\comdom_{\mathrm{i}}} \mu\mathbf{H}_{\mathrm{s}} \cdot \harmt_h\,\mathrm{d} \comdom
\end{align*}
for all $\mathbf{w}_h\in S_p^1(\comdom_{\mathrm{c}})$, $w_h\in S_p^0(\comdom)$ and $\harmt_h \in W_h(\comdom_{\mathrm{i}}) $ that must also satisfy the interface condition $\mathbf{w}_h \times \mathbf{n} = -\harmt_h \times \mathbf{n}$. As in the continuous case, it is necessary to add a gauging condition to ensure the uniqueness of the vector field $\mathbf{T}_h$. In our tests we will use tree-cotree gauging \cite{Albanese_1988aa}, that was adapted to the IGA setting in \cite{Kapidani_2022aa}.

What remains to be done is to construct a suitable discrete space $W_h(\comdom_{\mathrm{i}})$, which is the main contribution of this paper, and the subject of \cref{sec:cohomology}.

\begin{remark}
It is important to remark that, imposing the interface conditions strongly (as in \cref{eq:interface_Hphi_disc}), is in fact equivalent to use as trial and test functions a triple in the constrained space
\begin{equation*}
\left \{(\mathbf{w}_h, w_h,\harmt_h) \in S_p^1(\comdom_{\mathrm{c}}) \times S_p^0(\comdom_{\mathrm{i}}) \times W_h(\comdom_{\mathrm{i}}) : \mathbf{w}_h \times \mathbf{n} = (\grad \, w_h + \harmt_h) \times \mathbf{n} \text{ on } \Gamma
\right\}.
\end{equation*}
In practice, one chooses as basis functions for $\mathbf{w}_h$ only those associated to control edges in the internal part of the conductor, without the interface; then, to respect the interface condition, the functions $\phi_h$ and $w_h$ (and similarly $\harm_h$ and $\harmt_h$) on the interface will also give a contribution to the conductor part of the domain. These contributions correspond to the non-diagonal blocks appearing in the matrix in \cite{Webb_1993aa}. \end{remark}

\section{Cohomology with splines}\label{sec:cohomology}

Classic cohomology theory developed under the intuition that one could make statements about the topology of a domain based on studying classes of vector fields living on the domain. Computationally, it is nevertheless intuitive that once a mesh is established, computing combinatorially on the mesh graph is much more convenient than performing numerical orthogonalization within the subspace of curl-free vector fields to obtain the cohomology group generators, i.e., a basis of $W_h(\comdom_{\mathrm{i}})$ in our setting.

All the main published approaches towards computing cohomology groups generators indeed rely on the existence of a duality pairing between geometric entities in the mesh and discrete field unknowns, and thus work for simplicial meshes within either edge elements (also known as Whitney forms), the cell method \cite{Tonti_2001aa}, or FIT \cite{Clemens_2001aa} approaches, for which constructing such duality pairings is quite straightforward. 
For high order FEM, cohomology generators are left unchanged with respect to the one computed with edge elements, which are contained in the space of high order curl-conforming finite elements. Employing a hierarchical basis, e.g. the one from \cite{Schoberl_2005aa,Lee_2003aa}, where the basis of the lowest order elements is a subset of the high order basis then greatly simplifies the computations. 
Unfortunately, this approach cannot be straightforwardly applied for high continuity splines, because they do not contain the low order (and low continuity) spaces. 
On the other hand, in this section we will show that fortunately one can directly apply the same algorithms used for simplicial meshes in \cite{dlotkoLazyCohomologyGenerators2014,Dlotko_2017aa,Dlotko_2019aa} to the hexahedral control mesh generated by high-order splines, and thus compute cohomology generators as a linear combination of relatively few basis functions in the space of curl-conforming splines, for arbitrary order.

\subsection{From simplicial to cubical homology}
We start with a proposition showing that the cohomology properties for splines are the same as those of low order finite elements in the control mesh, the proof being a trivial application of the commutative isomorphisms in \eqref{eq:isomorphisms}. As a consequence, the cohomology computations for splines can be performed directly on the control mesh.

\begin{proposition}\label{prop:isomorphisms}
	Let $\tilde{H}(\comdom) = \mathrm{Ker}(\curl) / \mathrm{Im}(\grad)$ and $\tilde{H}(\comdom_C) = \mathrm{Ker}(\curl) / \mathrm{Im}(\grad)$ the first cohomology groups for splines in $\comdom$ and for low order finite elements in $\comdom_C$, respectively obtained from the top and bottom sequence in \eqref{eq:isomorphisms}. Then it holds that $ \dim(\tilde{H}(\comdom)) = \dim (\tilde{H}(\comdom_C)) = N_H$. Moreover, the functions $\{\harmt_{h,j}\}_{j=1}^{N_H} \subset S_p^1(\comdom)$ are generators of the first cohomology group $\tilde{H}(\comdom)$ if and only if $\{I_h^1(\harmt_{h,j})\}_{j=1}^{N_H} \subset Z_h^1(\comdom_C)$ are generators of the first cohomology group $\tilde{H}(\comdom_C)$.
\end{proposition}

After this result, the main critical issue in translating the algorithms from the literature to the setting of high order spline spaces is the fact that the control mesh we need to work with in the multipatch setting is an unstructured hexahedral one, whereas all the main algebraic tools are well established for simplicial homology. Luckily, the reasons for this are practical rather than theoretical, as laid out in the classic textbook by Hatcher~\cite{hatcherAlgebraicTopology2001}:

\vspace{12pt}

\begin{quotation}
    ``For a cell complex $X$ one has chain groups $C_n(X)$ which are free abelian groups with basis the n cells of $X$, and there are boundary homomorphisms $\partial_n :C_n(X)\mapsto C_{n-1}(X)$, in terms of which one defines the homology group $H_n(X) = \mathrm{Ker}(\partial_n) / \mathrm{Im}(\partial_{n+1})$. The major difficulty is how to define $\partial_n$ in general. [...]
The best solution to this problem seems to be to adopt an indirect approach.
Arbitrary polyhedra can always be subdivided into special polyhedra called simplices [...].
For simplices there is no difficulty in defining boundary maps or in handling orientations. So one obtains a homology theory, called simplicial homology, for cell complexes built from simplices.''
\end{quotation}

\vspace{12pt}
\noindent Nevertheless, in our case moving back and forth between an hexahedral control mesh and its simplicial subdivision would be very cumbersome. We wish instead to compute the cohomology generators directly on the control mesh. The theory behind the equivalence of the two approaches has been treated in the algebraic topology literature before, under the name of cubical homology (see e.g. \cite[Chap. 2]{kaczynski_computational_2004}).

The main hurdle with the boundary operators in this setting is achieving consistent global orientation of cells:
once orientation of all $n$-cells is established, the boundary operators $\partial_k : C_k(\comdom) \mapsto C_{k-1}(\comdom)$ on the whole complex can be represented by an incidence matrix $\mathbb{B}^k$ for which the entry with indices $i$, $j$ is formally given by:
\begin{equation}
  \mathbb{B}^k_{ij} = 
   \begin{cases}
    +1 & \text{if  the $i$-th ($k-1$)-cell is in the boundary of} \\ 
          & \text{the $j$-th $k$-cell with positive orientation;} \\
     -1 & \text{if  the $i$-th ($k-1$)-cell is in the boundary of} \\ 
         & \text{the $j$-th $k$-cell with negative orientation;} \\
      0 & \text{otherwise.}
   \end{cases}
 \label{eq:b-operator}
 \end{equation}

To make sense of this formal definition and solve the problem of constructing all boundary operators $\{\partial_k\}_{k=1}^{3}$ in our setting, one can proceed by using the local (patchwise structured) mesh topology information combined with the global ordering of vertices, by enumeration in any meshing output, as shown in Algorithm~\ref{algo:boundary-operators}, which is applied to each element of each patch in the control mesh (with the slight abuse of nomenclature remarked in \cref{sec:iga} when generalizing to multipatch.)
\begin{algorithm}\caption{Construction of boundary operators}
    \label{algo:boundary-operators}
\begin{enumerate}[1)]
    \item Each edge (or $1$-cell) is an ordered list of vertices $[v_1,v_2]$. $\mathbb{B}^1$ is then trivially constructed as in the simplicial mesh case.
    \item Each quadrilateral (or $2$-cell) is an ordered list of vertices $[v_1,v_2,v_3,v_4]$. There are two edges with $v_1$ in their boundary. \emph{We label these two edges as $e_1$ and $e_2$, according to which one has the lowest subscript in the second vertex in the pair $[v_1,v_\star]$ (where $\star\in\{2,3,4\}$) and always choose them to induce the outer orientation of the $2$-cell via $e_1\wedge e_2$}. The outer orientation then induces the inner orientation and then entries of $\mathbb{B}^2$ can be easily computed.
    \item For each hexahedron (or $3$-cell) one can then define the positive outer orientation via the outer normal vector field, and easily compute if the outer orientations of the six facets in its boundary either coincide or are opposite to the outer unit normal.
\end{enumerate}
\end{algorithm}
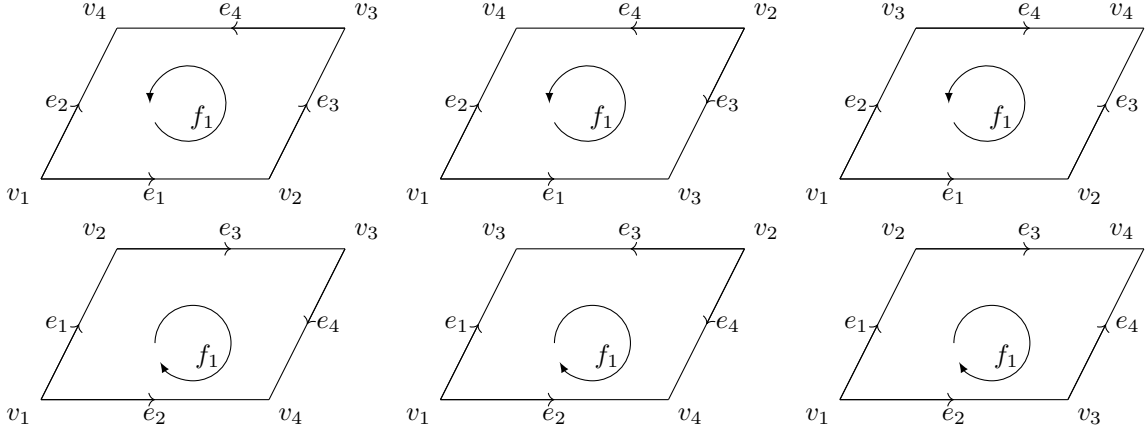
\begin{figure}\centering
\scalebox{1.0}{
\begin{tikzpicture}
\coordinate (v1) at (0,0);
\coordinate (v2) at (3,0);
\coordinate (v3) at (4,2);
\coordinate (v4) at (1,2);

\draw (v1) -- node[below] {$e_1$} (v2);
\draw (v1) -- node[left] {$e_2$} (v4);
\draw (v2) -- node[right] {$e_3$} (v3);
\draw (v3) -- node[above] {$e_4$} (v4);

\node[below left] at (v1) {$v_1$};
\node[below right] at (v2) {$v_2$};
\node[above right] at (v3) {$v_3$};
\node[above left] at (v4) {$v_4$};

\draw[->] (v1) -- node[below] {} ++($0.5*(v1)+0.5*(v2)$);
\draw[->] (v1) -- node[left] {} ++($0.5*(v1)+0.5*(v4)$);
\draw[->] (v2) -- node[right] {} ++($0.5*(v3)-0.5*(v2)$);
\draw[->] (v3) -- node[above] {} ++($0.5*(v4)-0.5*(v3)$);
\draw[-latex,color=black] (1.5,0.75) arc (-150:180:5mm) node[near start,above] {$f_1$};
\end{tikzpicture}
\begin{tikzpicture}
\coordinate (v1) at (0,0);
\coordinate (v2) at (4,2);
\coordinate (v3) at (3,0);
\coordinate (v4) at (1,2);

\draw (v1) -- node[below] {$e_1$} (v3);
\draw (v1) -- node[left] {$e_2$} (v4);
\draw (v2) -- node[right] {$e_3$} (v3);
\draw (v2) -- node[above] {$e_4$} (v4);

\node[below left] at (v1) {$v_1$};
\node[above right] at (v2) {$v_2$};
\node[below right] at (v3) {$v_3$};
\node[above left] at (v4) {$v_4$};

\draw[->] (v1) -- node[below] {} ++($0.5*(v1)+0.5*(v3)$);
\draw[->] (v1) -- node[left] {} ++($0.5*(v1)+0.5*(v4)$);
\draw[->] (v2) -- node[right] {} ++($0.5*(v4)-0.5*(v2)$);
\draw[->] (v2) -- node[above] {} ++($-0.5*(v2)+0.5*(v3)$);
\draw[-latex,color=black] (1.5,0.75) arc (-150:180:5mm) node[near start,above] {$f_1$};
\end{tikzpicture}
\begin{tikzpicture}
\coordinate (v1) at (0,0);
\coordinate (v4) at (4,2);
\coordinate (v2) at (3,0);
\coordinate (v3) at (1,2);

\draw (v1) -- node[below] {$e_1$} (v2);
\draw (v1) -- node[left] {$e_2$} (v3);
\draw (v2) -- node[right] {$e_3$} (v4);
\draw (v3) -- node[above] {$e_4$} (v4);

\node[below left] at (v1) {$v_1$};
\node[below right] at (v2) {$v_2$};
\node[above left] at (v3) {$v_3$};
\node[above left] at (v4) {$v_4$};

\draw[->] (v1) -- node[below] {} ++($0.5*(v1)+0.5*(v2)$);
\draw[->] (v1) -- node[left] {} ++($0.5*(v1)+0.5*(v3)$);
\draw[->] (v2) -- node[right] {} ++($0.5*(v4)-0.5*(v2)$);
\draw[->] (v3) -- node[above] {} ++($0.5*(v4)-0.5*(v3)$);
\draw[-latex,color=black] (1.5,0.75) arc (-150:180:5mm) node[near start,above] {$f_1$};
\end{tikzpicture}
}

\scalebox{1.0}{
    \begin{tikzpicture}
\coordinate (v1) at (0,0);
    \coordinate (v4) at (3,0);
    \coordinate (v3) at (4,2);
    \coordinate (v2) at (1,2);
    
\draw (v1) -- node[left] {$e_1$} (v2);
    \draw (v1) -- node[below] {$e_2$} (v4);
    \draw (v4) -- node[right] {$e_4$} (v3);
    \draw (v3) -- node[above] {$e_3$} (v2);
    
\node[below left] at (v1) {$v_1$};
    \node[below right] at (v4) {$v_4$};
    \node[above right] at (v3) {$v_3$};
    \node[above left] at (v2) {$v_2$};
    
\draw[->] (v1) -- node[below] {} ++($0.5*(v1)+0.5*(v4)$);
    \draw[->] (v1) -- node[left] {} ++($0.5*(v1)+0.5*(v2)$);
    \draw[->] (v3) -- node[right] {} ++($-0.5*(v3)+0.5*(v4)$);
    \draw[->] (v2) -- node[above] {} ++($-0.5*(v2)+0.5*(v3)$);
    \draw[-latex,color=black] (1.5,0.75) arc (180:-150:5mm) node[near end,above] {$f_1$};
    \end{tikzpicture}
    \begin{tikzpicture}
\coordinate (v1) at (0,0);
        \coordinate (v2) at (4,2);
        \coordinate (v4) at (3,0);
        \coordinate (v3) at (1,2);
        
\draw (v1) -- node[left] {$e_1$} (v3);
        \draw (v1) -- node[below] {$e_2$} (v4);
        \draw (v2) -- node[right] {$e_4$} (v4);
        \draw (v2) -- node[above] {$e_3$} (v3);
        
\node[below left] at (v1) {$v_1$};
        \node[above right] at (v2) {$v_2$};
        \node[below right] at (v4) {$v_4$};
        \node[above left] at (v3) {$v_3$};
        
\draw[->] (v1) -- node[below] {} ++($0.5*(v1)+0.5*(v4)$);
        \draw[->] (v1) -- node[left] {} ++($0.5*(v1)+0.5*(v3)$);
        \draw[->] (v2) -- node[right] {} ++($0.5*(v3)-0.5*(v2)$);
        \draw[->] (v2) -- node[above] {} ++($-0.5*(v2)+0.5*(v4)$);
        \draw[-latex,color=black] (1.5,0.75) arc (180:-150:5mm) node[near end,above] {$f_1$};
    \end{tikzpicture}
    \begin{tikzpicture}
\coordinate (v1) at (0,0);
    \coordinate (v4) at (4,2);
    \coordinate (v3) at (3,0);
    \coordinate (v2) at (1,2);
    
\draw (v1) -- node[left] {$e_1$} (v2);
    \draw (v1) -- node[below] {$e_2$} (v3);
    \draw (v2) -- node[above] {$e_3$} (v4);
    \draw (v3) -- node[right] {$e_4$} (v4);
    
\node[below left] at (v1) {$v_1$};
    \node[below right] at (v3) {$v_3$};
    \node[above left] at (v2) {$v_2$};
    \node[above left] at (v4) {$v_4$};
    
\draw[->] (v1) -- node[below] {} ++($0.5*(v1)+0.5*(v3)$);
    \draw[->] (v1) -- node[left] {} ++($0.5*(v1)+0.5*(v2)$);
    \draw[->] (v3) -- node[right] {} ++($0.5*(v4)-0.5*(v3)$);
    \draw[->] (v2) -- node[above] {} ++($0.5*(v4)-0.5*(v2)$);
    \draw[-latex,color=black] (1.5,0.75) arc (180:-150:5mm) node[near end,above] {$f_1$};
    \end{tikzpicture}
}
\caption{The possible configurations for a quadrilateral according to the global ordering of its vertices, up to cyclic permutation, together with orientation and induced ordering of edges $[e_1,e_2,e_3,e_4]$ and the orientation of the quadrilateral itself as computed by our Algorithm~\ref{algo:boundary-operators}.}
\label{fig:quad_orientation}
\end{figure}

For the boundary operator $\partial_2$ we note that the procedure above can be cumbersome to understand, so we elucidate it further by giving all the possible configurations with arrows indicating orientations in Figure~\ref{fig:quad_orientation}.
Once the boundary operators are constructed we state the following proposition, with trivial proofs, which we will use in the following.

\begin{proposition}\label{lemmino}
    \textbf{Co-boundary operators:} The transpose matrices given by \eqref{eq:b-operator} for the control mesh $V_C$, i.e. $$(\mathbb{B}_1)^\top\;(\mathbb{B}_2)^\top,\; (\mathbb{B}_3)^\top,$$ respectively encode the operators $$\grad: Z_h^0(V_C) \mapsto Z_h^1(V_C),\;\curl: Z_h^1(V_C) \mapsto Z_h^2(V_C),\; \div: Z_h^2(V_C) \mapsto Z_h^3(V_C),$$ in diagram \eqref{eq:isomorphisms}.
\end{proposition}

To see why \cref{lemmino} works, let us consider matrix $(\mathbb{B}_1)^\top$: it suffices to note that its restriction to each patch $V^{(j)}$ coincides with the matrix $\mathbb{G}$ up to reordering of basis functions and sign changes, depending on the orientation. The same applies to the restriction of $(\mathbb{B}_2)^\top$ to any $V^{(j)}$ and $\mathbb{C}$, and to the restriction of $(\mathbb{B}_3)^\top$ to any $V^{(j)}$ and $\mathbb{D}$, through tensorization arguments. Furthermore, clearly $\mathbb{B}_{k+1} \mathbb{B}_k = 0$ (i.e. an identically null matrix, see \cite[Prop. 2.37]{kaczynski_computational_2004}).

\subsection{The cohomology algorithm}
We recall that the goal is to compute a basis of the space $W_h(\comdom_{\mathrm{i}})$, as introduced in Section~\ref{subsec:weak_IGA}. Before explaining the algorithm to compute this basis using the control mesh, we introduce some necessary notation. We will assume that each patch is contained either in the conductor or in the insulator subdomain, and therefore the interface $\Gamma$ is given as the union of interfaces between patches. We will keep the notation $\comdom_C$ for the control mesh of the whole domain, and denote by $\comdom_{\mathrm{c},C}$ and $\comdom_{\mathrm{i},C}$ the control mesh for each subdomain, and by $\Gamma_C$ the control mesh defining the interface, which is in fact the interface between $\comdom_{\mathrm{c},C}$ and $\comdom_{\mathrm{i},C}$.

Building on the insight of \cref{prop:isomorphisms} and \cref{lemmino}, a basis of $W_h(\comdom_{\mathrm{i}})$ can be obtained from a set of generators of the cohomology group $\tilde{H}(V_{\mathrm{i},C})$ computed on the control mesh. The computation of these generators is a straightforward application of the method developed in \cite{dlotko_lean_2018}. We will not dive too deep into detail within the algorithmic steps, as they are an adaption from simplicial to cubical homology, and refer to that paper for the details.

We first remark that the dimension of the cohomology group is not known a priori. This is in fact induced by the cohomology on $\Gamma$ (or $\Gamma_C$), i.e., the space generated as the quotient $$\widetilde{H}(\Gamma_C) = \mathrm{Ker}(\mathrm{curl}_\Gamma) / \mathrm{Im}(\grad_\Gamma),$$ with reference to the differential operators used in \eqref{eq:spline_derham_gamma}, or analogous spaces for finite elements on $\Gamma_C$. A basis of cohomology generators on $\Gamma_C$ is computed with algorithms developed in \cite{kapidaniComputationRelative1Cohomology2016}, itself based on \cite{hiptmairGeneratorsGammah2002}. Then, each of these generators is given as input to \cref{algo:generators}, to obtain the cohomology generators for $\tilde{H}(\comdom_{\mathrm{i},C})$.

\begin{algorithm}\caption{Construction of a cohomology generator on the control mesh}
	\begin{enumerate}[1)]
\item Starting from a cohomology generator on $\Gamma_C$, construct a discrete function $\tilde{\boldsymbol{j}}\in Z_h^2(\comdom_{\mathrm{c},C})$ such that $\tilde{\boldsymbol{j}}\cdot\mathbf{n}=0$ on $\Gamma$ and $\div(\tilde{\boldsymbol{j}}) = 0$ in $\comdom_{\mathrm{c},C}$. Extend $\tilde{\boldsymbol{j}}$ to $\comdom_{\mathrm{i},C}$ by zero.\label{algo:step1}
		\item Remove a spanning tree of edge degrees of freedom (DoFs) from $Z_h^1(V_C)$, identified through a breadth first search algorithm (BFS~\cite[Chap. 23.2]{Cormen_2001aa}). \label{algo:step2}
		\item Define $\mathbb{C}_\square$ as the row echelon form of the square matrix obtained by eliminating the corresponding columns of ${(\mathbb{B}^2)}^\top$ and solve the system $$\mathbb{C}_\square \tilde{\mathbf{h}} = \tilde{\mathbf{j}}$$ by backward substitution, where $\tilde{\mathbf{j}}$ are the DoFs of $\tilde{\boldsymbol{j}}$.\label{algo:step3}
		\item Compute $\tilde{\mathbf{h}}_\mathrm{i}$ the restriction of $\tilde{\mathbf{h}}$ to $\comdom_{\mathrm{i},C}$.\label{algo:output}
\end{enumerate}
	\label{algo:generators}
\end{algorithm}

The steps of \cref{algo:generators} can be summarised as follows: In step \ref{algo:step1} we compute a divergence free current inside the conductor for each cohomology generator we are seeking to construct. Then, in step~\ref{algo:step2} we remove a spanning tree from the whole control mesh, in such a way that the matrix $\mathbb{C}_\square$ computed in step~\ref{algo:step3} is square and non-singular. This procedure is akin to high-order versions of tree-cotree techniques used in the gauging of magnetostatics formulations in \cite{Merkel_2022aa,Los-Santos_2023aa}.
The solution of the linear system at this step is a function defined in the whole domain $\comdom_C$ and with nonzero curl equal to $\tilde{\boldsymbol{j}}$, thus it cannot be a gradient. In step~\ref{algo:output} we restrict this function to $\comdom_{\mathrm{i},C}$. From the assumptions on $\tilde{\boldsymbol{j}}$ this restriction has zero curl, but from the previous step it cannot be a gradient, and therefore it belongs to $\tilde{H}(\comdom_{\mathrm{i},C})$. Finally, applying the inverse of the isomorphism $I_h^1$
in \eqref{eq:isomorphisms} to the computed generators, we obtain the generators for $\tilde{H}(\comdom_{\mathrm{i}})$, and therefore a basis for $W_h(\comdom_{\mathrm{i}})$. This operation is inexpensive, since the choice of the B-spline basis guarantees that the isomorphisms are given by diagonal matrices. 

It is important to note that the final cohomology generators will be curl-conforming spline functions in $W_h(\comdom_{\mathrm{i}}) \subset S_p^1(\comdom_{\mathrm{i}})$, but will still have by construction very local support within the mesh. We also remark that an expensive gaussian elimination procedure is never performed to obtain $\mathbb{C}_\square$ in steps~\ref{algo:step2}--\ref{algo:step3}, which go by the name of spanning tree technique (STT). Rather, the STT itself provides the recipe for the backwards substitution steps by exploiting graph locality in the tree construction algorithm. Moreover, although the algorithm is written for a single generator, in case there are many the same spanning tree in step~\ref{algo:step2} and the same system matrix $\mathbb{C}_\square$ in step~\ref{algo:step3} can be used for all of them, for which one has to solve a system with multiple right-hand sides. We will back this claim up in \cref{sec:numerics} through numerical experiments, which corroborate the lowest order results from \cite{Dlotko_2019aa}.

We close the section with a remark regarding the computational complexity of Algorithm~\ref{algo:generators}.

\begin{remark}
    In the assumption that $\mathbb{C}_\square$ is always available after elimination of tree DoFs, all steps in Algorithm~\ref{algo:generators} are at most linear in the number of unknowns in the $S_p^1(V)$ space. This assumption is usually valid, although some pathological cases exist, as discussed at length in \cite{dlotkoCriticalAnalysisSpanning2010}. Recent work in \cite{pitassiInvertingDiscreteCurl2022} strongly suggests (without proof) that local post-processing procedures to step \ref{algo:step2} can finally eliminate these corner cases. Conveniently, the algorithms therein introduced also easily translate to cubical cohomology, and to our high-order spline setting. Nevertheless, further discussion of the details goes beyond the scope of present article.

    Furthermore, the machinery outlined in Algorithm~\ref{algo:generators} will actually provide double the number of generators with respect to the one needed since $\dim{\tilde{H}(\Gamma_C)} = 2\dim{W_h(\comdom_{\mathrm{i},C})}$. This choice was made for the sake of simplicity in our exposition, and motivates the nomenclature \emph{lazy generators} which has been used, e.g., in \cite{dlotkoLazyCohomologyGenerators2014}. Details on how to mend this issue at negligible computational cost (i.e. without computing expensive double integrals for linking numbers between homology group generators as done in \cite{rodriguezConstructionFiniteElement2013}) can be found in \cite{Dlotko_2019aa}.
\end{remark}

\section{Numerical Results}\label{sec:numerics}
\definecolor{colorA}{RGB}{202,0,32}  
\definecolor{colorB}{RGB}{244,165,130}  
\definecolor{colorC}{RGB}{146,197,222}  
\definecolor{colorD}{RGB}{5,113,176}  
The framework constructed in the previous sections is hereafter supported by numerical tests, first on a toy problem of academic interest and then on a benchmark developed by experimentalists to test the accuracy of numerical methods for the eddy currents problem. Finally, we test the computational complexity of the algorithm in Section~\ref{sec:cohomology} with respect to the number of holes.
For the implementation of isogeometric analysis in Matlab in the following computations, we use the open source package GeoPDEs \cite{Vazquez_2016aa} and the NURBS toolbox \cite{NURBS_2021aa}. The cohomology generators are computed using the Topoprocessor toolbox \cite{Dlotko_2017aa,bkapidaniREADME2021}. All the aforementioned toolboxes are open source software available freely to researchers online. All computations are carried out in Matlab\textsuperscript{\textregistered} R2022b on a 6-core machine (Intel\textsuperscript{\textregistered}
 Core\texttrademark{} i7-5820K CPU) with \SI{16}{GB} RAM. 

\subsection{Hollow Cylinder}
\begin{figure}
    \center
\begin{tikzpicture}
            \node at (0,0) {\includegraphics[width=0.4\columnwidth,trim=2cm 4.5cm 7cm 6cm, clip]{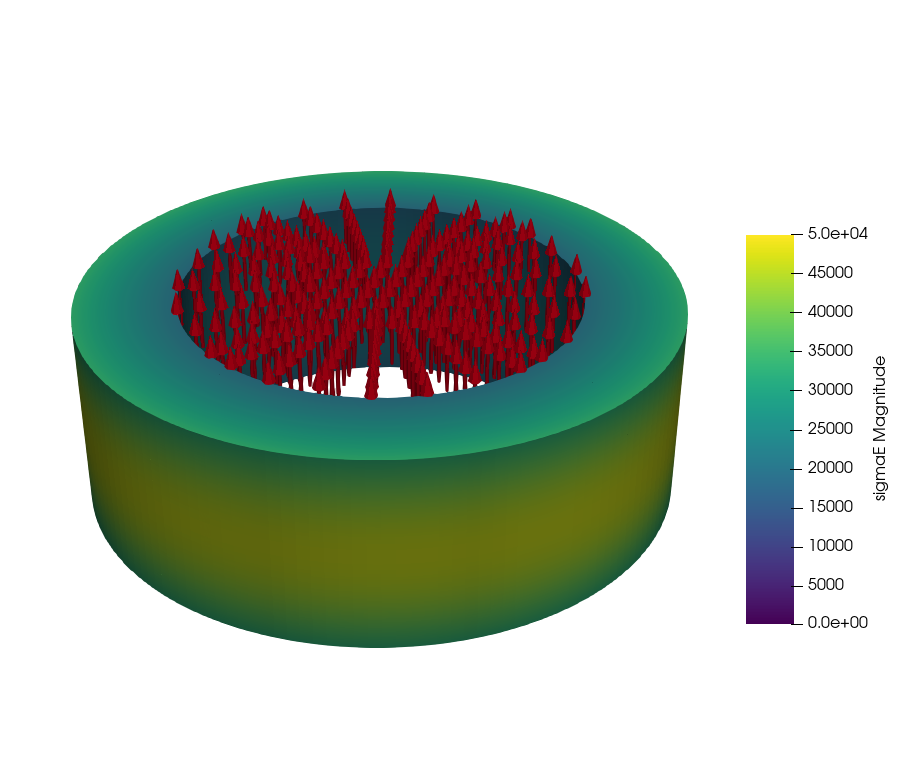}};
\end{tikzpicture}
    \caption{Support of the computed basis function in $W_h(\comdom_{\mathrm{i}})$, spanned by the red arrows.}
    \label{fig:cylinder_harmonic_fields}
\end{figure}

As a first toy problem to numerically test if our algorithms handle correctly curved geometries, we consider a conducting hollow cylinder with a conductivity of $\sigma=\SI{3.256e7}{\siemens\per\meter}$ which is placed concentrically inside a coil as visualized in Figure~\ref{fig:cylinder_geometry}. 
The geometric parameters of the cylinder and the coil are given by $R_{\mathrm{i}}=\SI{1}{cm}$, $R_{\mathrm{o}}=\SI{1.5}{cm}$, $R_{\mathrm{c,i}}=\SI{1.8}{cm}$, $R_{\mathrm{c,o}}=\SI{2}{cm}$ and $h=\SI{1}{cm}$. The amplitude of the total current in the coil (ampere-turn) is $\SI{20}{\ampere}$ in the angular direction at a frequency of $\omega = 2\pi\cdot \SI{50}{\hertz}$ and is assumed uniformly distributed across the cross-section of the coil. The conductor is described with 4 patches and discretized with $864$ elements, while the whole domain is given by $39$ patches and discretized with $9288$ elements. The problem is discretized with splines of degree $p=3$, which gives a total of 26792 unknowns for the magnetic $\vecs{H}-\phi$ and the $\vecs{T}-\Omega$ formulations, and $72520$ unknowns for the electric $\vecs{A}-\varphi$-formulation. The resulting eddy currents in the conductor are shown in Figure~\ref{fig:cylinder_sigmaE} for the $\vecs{H}-\phi$ and the $\vecs{T}-\Omega$ formulations, and as a comparison for the well-known (but more computationally expensive) $\vecs{A}-\varphi$-formulation where the homology of the conductor requires no special treatment. We observe that the two magnetic formulations yield the same result up to round-off errors, while the electric $\vecs{A}-\varphi$ formulation gives very similar results. 
 Since there is a single hole, there is only one basis function in $W_h(\comdom_{\mathrm{i}})$. 
The number $\Nharm$ of basis functions of $S_p^1(\comdom_{\mathrm{i}})$ which are linearly combined to obtain this function, and the number $\Ngrad$ of basis functions of $S_p^0(\comdom_{\mathrm{i}})$, are $\Nharm=337$ and $\Ngrad=24457$.
The support of the only basis function in $W_h(\comdom_{\mathrm{i}})$ is shown in Figure~\ref{fig:cylinder_harmonic_fields} through a vector field plot, and it is very localized inside the hole in the conductor. 
\begin{figure}
    \center
    \usetikzlibrary{shapes.geometric}

\begin{tikzpicture}[]
\def\ri{1}
\def\ro{1.5}
\def\h{1}

\def\ric{2.2}
\def\roc{2.5}
\def\hc{1}

{\color{purple}
\draw [thick](-\ro,0) -- (-\ro,\h);
\draw [thick](\ro,0) -- (\ro,\h);
\draw [thick](-\ro,0) arc (180:360:\ro\space and \ro/3);          \draw[thick,dashed] (\ro,0) arc (-\ro:180:\ro\space and \ro/3);  \draw [thick](+\ro,+\h) arc (-\ro:360:\ro\space and \ro/3) -- cycle;         \draw [thick](-\ri,0) -- (-\ri,\h);
\draw [thick](\ri,0) -- (\ri,\h);
\draw [thick](-\ri,0) arc (180:360:\ri\space and \ri/3);
\draw[thick,dashed] (1,0) arc (-\ri:180:\ri\space and \ri/3);
\draw [thick](+\ri,+\h) arc (-\ri:360:\ri\space and \ri/3) -- cycle;
\draw [thick,<->] (0,0) -- node[fill=white,inner sep=-0.5pt] {\scriptsize$R_{\mathrm{i}}$} (\ri,0);
\draw [thick,<->,yshift=-5pt] (0,-\roc/3) -- node[fill=white,scale=0.9,inner sep=1pt] {\scriptsize$R_{\mathrm{o}}$} (\ro,-\roc/3);
}

{\color{gray}
\draw [thick](-\roc,0) -- (-\roc,\hc);
\draw [thick](\roc,0) -- (\roc,\hc);
\draw [thick](-\roc,0) arc (180:360:\roc\space and \roc/3);          \draw[thick,dashed] (\roc,0) arc (-\roc:180:\roc\space and \roc/3);  \draw [thick](+\roc,+\hc) arc (-\roc:360:\roc\space and \roc/3) -- cycle;         \draw [thick](-\ric,0) -- (-\ric,\hc);
\draw [thick](\ric,0) -- (\ric,\hc);
\draw [thick](-\ric,0) arc (180:360:\ric\space and \ric/3);
\draw[thick,dashed] (\ric,0) arc (-\ric:180:\ric\space and \ric/3);
\draw [thick](+\ric,+\hc) arc (-\ric:360:\ric\space and \ric/3) -- cycle;
\draw [thick,<->] (0,0) -- node[fill=white,inner sep=-0.5pt] {\scriptsize$R_{\mathrm{c,i}}$} (-\ric,0);
\draw [thick,<->,yshift=-5pt] (0,-\roc/3) -- node[fill=white,scale=0.9,inner sep=1pt] {\scriptsize$R_{\mathrm{c,o}}$} (-\roc,-\roc/3);
\draw [thick,<->,xshift=5pt] (\roc,\hc) -- node[fill=white,scale=0.9, inner sep=1pt,anchor=west] {$h$} (\roc,0);
}
\end{tikzpicture}
     \caption{Sketch of the setup where a conducting hollow cylinder (red) with $\sigma=\SI{3.256e7}{\siemens\per\meter}$ is placed inside a concentric coil (gray).}
    \label{fig:cylinder_geometry}
\end{figure}
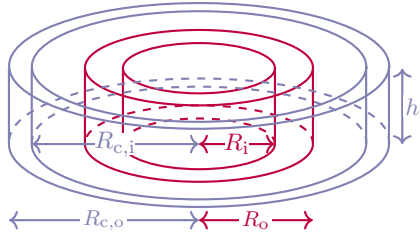 
\begin{figure}
    \center
        \begin{subfigure}[c]{0.25\columnwidth}
        \center
        \begin{tikzpicture}
            \node at (0,0) {\includegraphics[height=3cm,trim=2cm 4.5cm 7cm 6cm, clip]{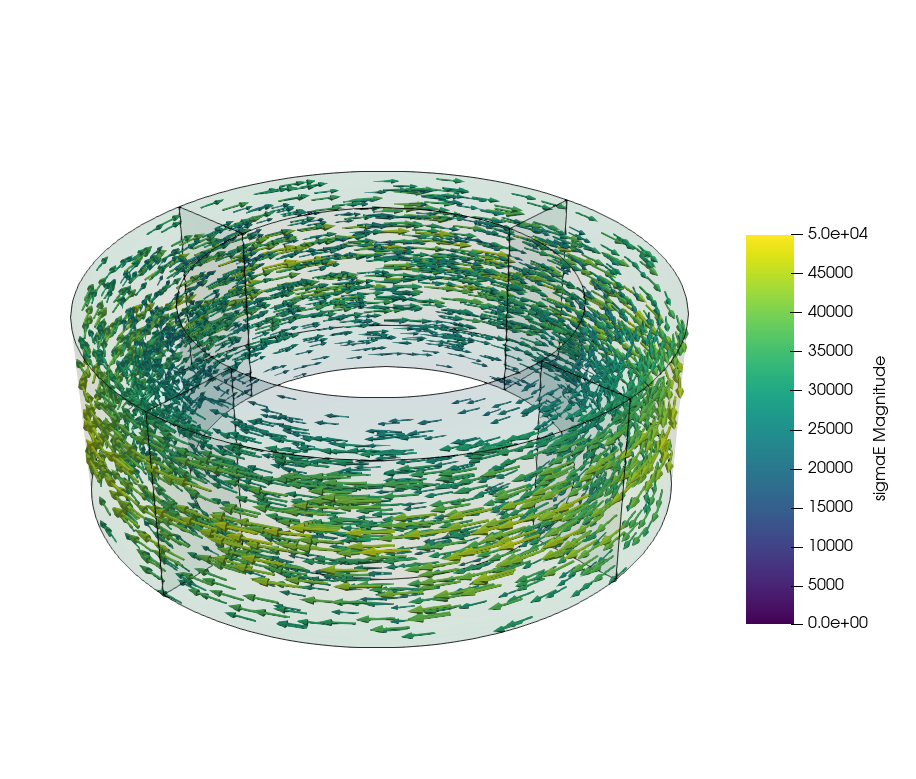}};
            \node [anchor=west] at (0,1.7) {\phantom{$\lvert\vecs{J}_\textrm{c}\rvert$ (\si{A\per\meter\squared})}};
        \end{tikzpicture}
        \\
        \begin{tikzpicture}
            \node at (0,0) {\includegraphics[height=3cm,trim=2cm 4.5cm 7cm 6cm, clip]{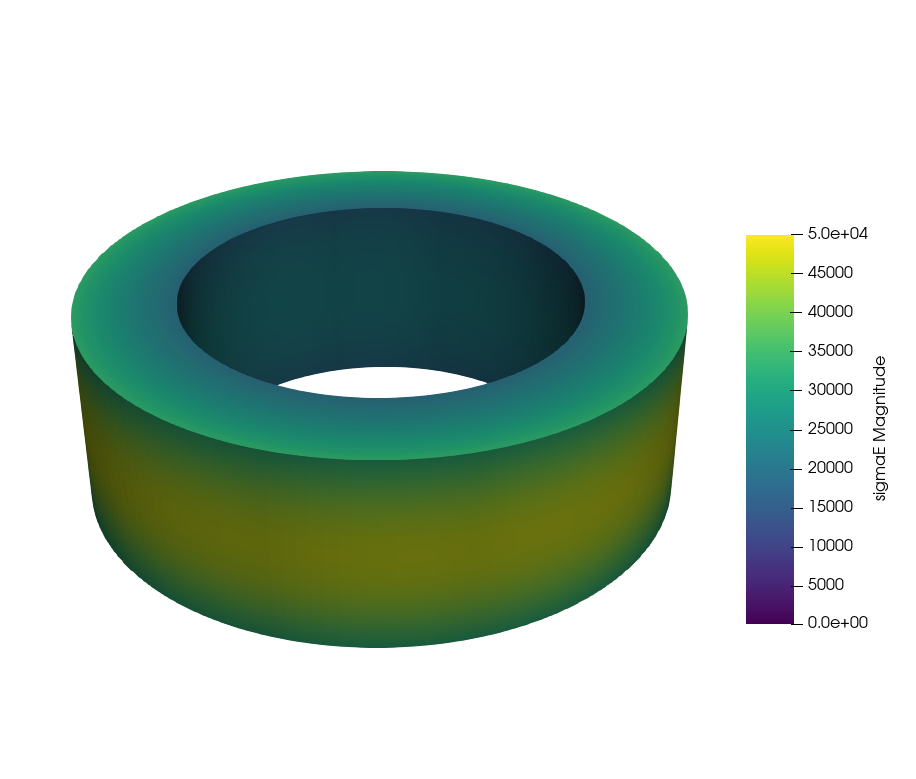}};
            \node [anchor=west] at (0,1.7) {\phantom{$\lvert\vecs{J}_\textrm{c}\rvert$ (\si{A\per\meter\squared})}};
        \end{tikzpicture}
        \subcaption{$\vecs{H}-\phi$}\label{fig:team7_J_Hphi}
    \end{subfigure}
    \hfill
    \begin{subfigure}[c]{0.25\columnwidth}
        \center
        \begin{tikzpicture}
            \node at (0,0) {\includegraphics[height=3cm,trim=2cm 4.5cm 7cm 6cm, clip]{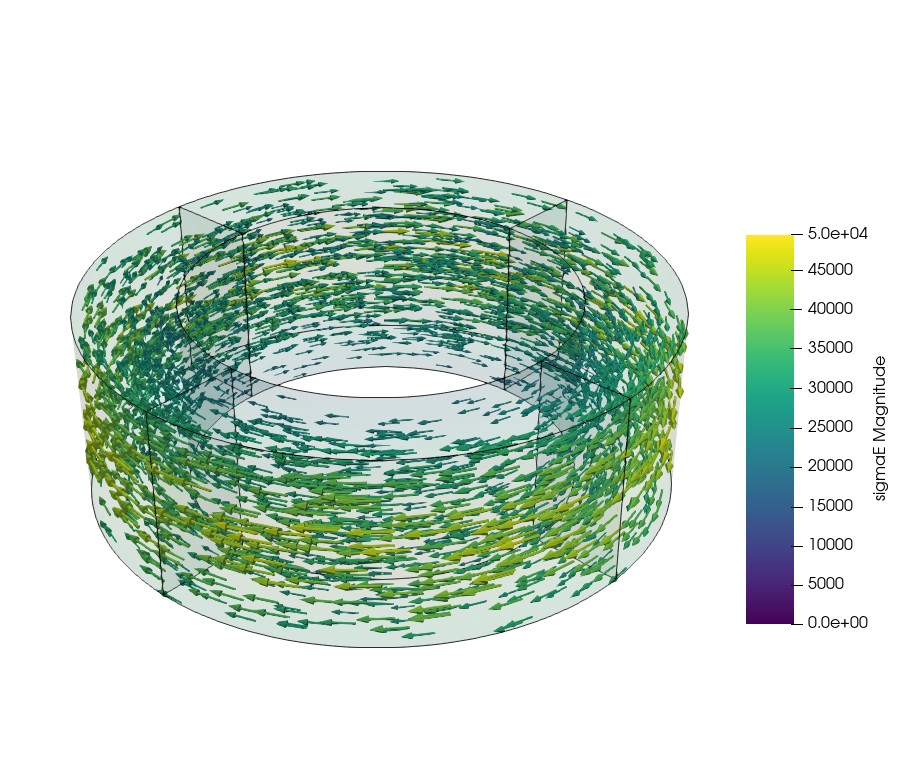}};
             \node [anchor=west] at (0,1.7) {\phantom{$\lvert\vecs{J}_\textrm{c}\rvert$ (\si{A\per\meter\squared})}};
        \end{tikzpicture}
        \\
        \begin{tikzpicture}
            \node at (0,0) {\includegraphics[height=3cm,trim=2cm 4.5cm 7cm 6cm, clip]{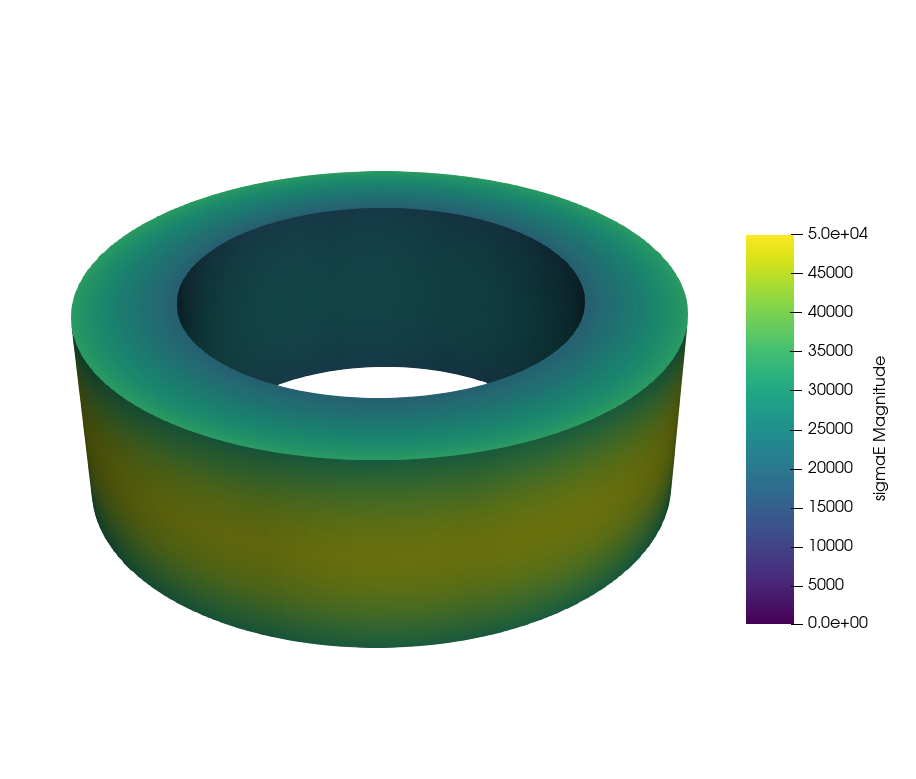}};
             \node [anchor=west] at (0,1.7) {\phantom{$\lvert\vecs{J}_\textrm{c}\rvert$ (\si{A\per\meter\squared})}};
        \end{tikzpicture}
        \subcaption{$\vecs{T}-\Omega$}
    \end{subfigure}
    \hfill
    \begin{subfigure}[c]{0.4\columnwidth}
        \center
        \begin{tikzpicture}
            \node at (0,0) {\includegraphics[height=3cm,trim=2cm 4.5cm 5cm 6cm, clip]{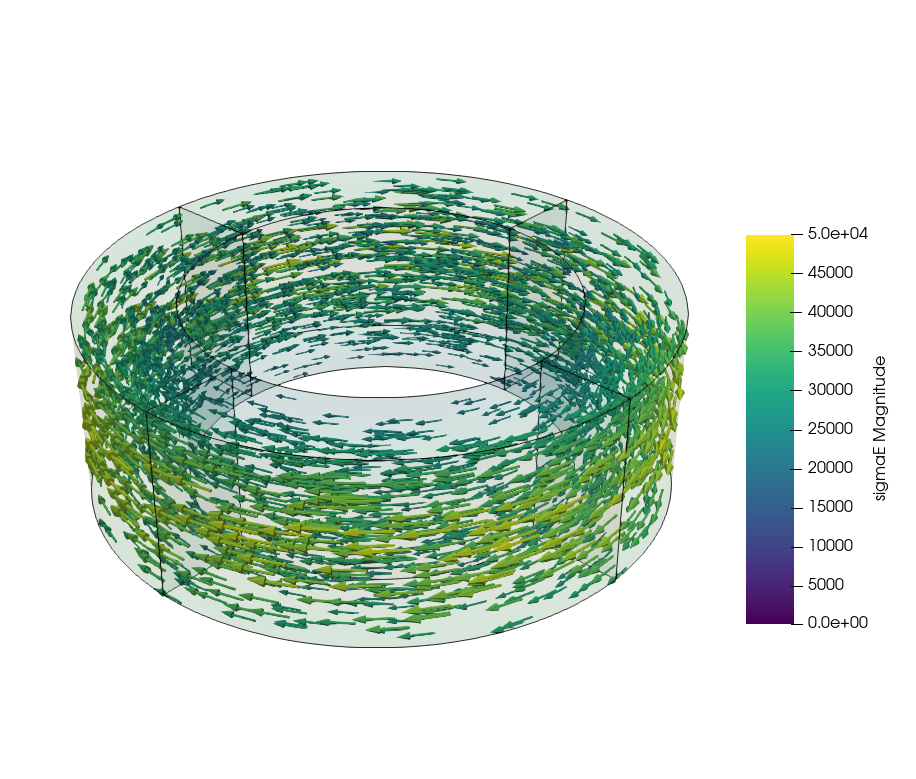}};
            \foreach \i in {0,1,...,5.01}{
                \draw (2.16,2/5*1.2*\i-1.32)--++(5pt,0) node [pos=1,anchor=west] {\small$\pgfmathprintnumber[fixed,precision=1]{\i}\cdot 10^4$};
             }
             \node [anchor=west] at (1.9,1.7) {$\lvert\vecs{J}_\textrm{c}\rvert$ (\si{A\per\meter\squared})};
        \end{tikzpicture}
        \\
        \begin{tikzpicture}
            \node at (0,0) {\includegraphics[height=3cm,trim=2cm 4.5cm 5cm 6cm, clip]{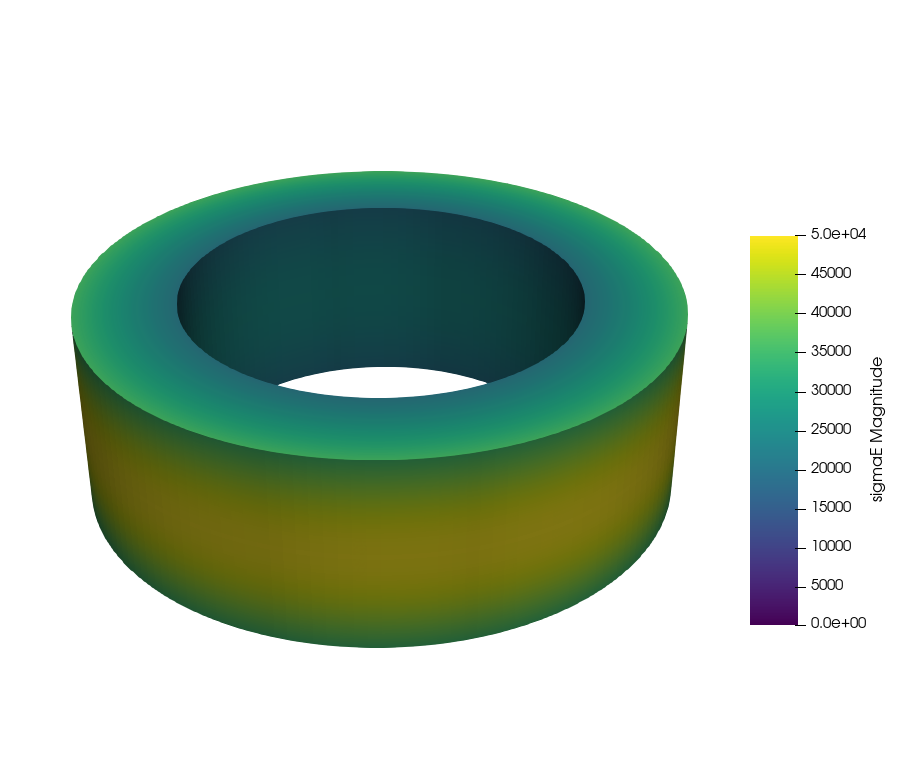}};
            \foreach \i in {0,1,...,5.01}{
                \draw (2.16,2/5*1.2*\i-1.32)--++(5pt,0) node [pos=1,anchor=west] {\small$\pgfmathprintnumber[fixed,precision=1]{\i}\cdot 10^4$};
             }
             \node [anchor=west] at (1.9,1.7) {$\lvert\vecs{J}_\textrm{c}\rvert$ (\si{A\per\meter\squared})};
        \end{tikzpicture}
        \subcaption{$\vecs{A}-\varphi$}
    \end{subfigure}
    \caption{Eddy currents $\vecs{J}_\mathrm{c}$ in the conducting hollow cylinder at $t = \tfrac{1}{4\omega}$ computed using the formulations in \cref{sec:maxwell} with basis functions of degree $p=3$. \textit{Top:} Direction of the currents $\vecs{J}_\textrm{c}$. \textit{Bottom:} Magnitude $\lvert\vecs{J}_\textrm{c}\rvert$ of the eddy currents.}
    \label{fig:cylinder_sigmaE}
\end{figure} 

\subsection{Team 7 Benchmark}
As a second example, we consider the TEAM benchmark problem 7 \cite{Fujiwara_1990aa}, where a conducting plate with a hole is placed in a non-uniform magnetic field, generated by a coil above the plate with a sinusoidal (total) current of the frequency $\omega=2\pi\cdot \SI{50}{\hertz}$ and the maximum amplitude $\SI{2742}{\ampere}$ at time $t=\SI{0}{\second}$. For details of the geometric setup we refer to \cite{Fujiwara_1990aa}.

\begin{figure}
    \center
    \begin{subfigure}[c]{0.27\columnwidth}
        \center
        \begin{tikzpicture}
            \node at (0,0) {\includegraphics[height=2.5cm,trim=0cm 5cm 7cm 9cm, clip]{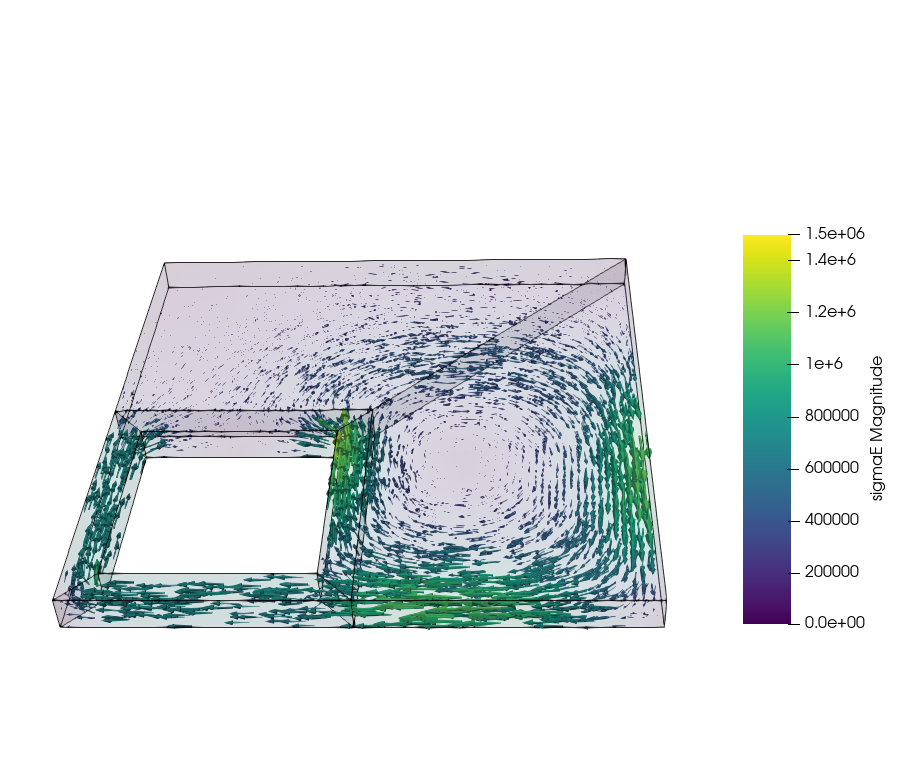}};
            \node [anchor=west] at (0,1.7) {\phantom{$\lvert\vecs{J}_\textrm{c}\rvert$ (\si{A\per\meter\squared})}};
        \end{tikzpicture}
        \\
        \begin{tikzpicture}
            \node at (0,0) {\includegraphics[height=2.5cm,trim=0cm 5cm 7cm 9cm, clip]{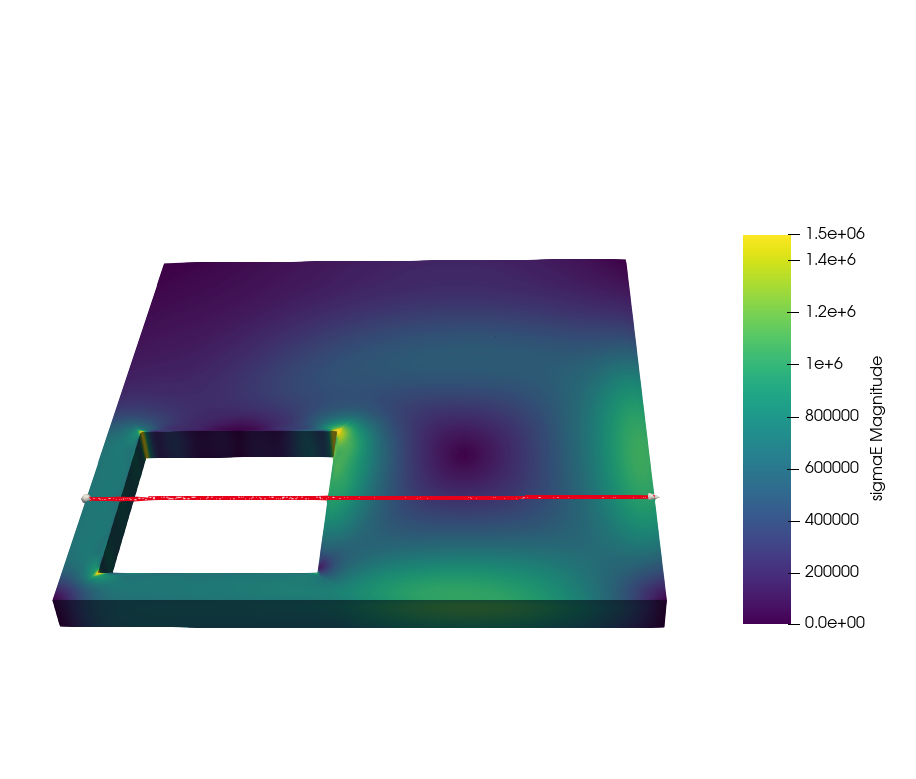}};
            \node [anchor=west] at (0,1.7) {\phantom{$\lvert\vecs{J}_\textrm{c}\rvert$ (\si{A\per\meter\squared})}};
        \end{tikzpicture}
        \subcaption{$\vecs{H}-\phi$}\label{fig:team7_J_Hphi}
    \end{subfigure}
    \hfill
    \begin{subfigure}[c]{0.27\columnwidth}
        \center
        \begin{tikzpicture}
            \node at (0,0) {\includegraphics[height=2.5cm,trim=0cm 5cm 7cm 9cm, clip]{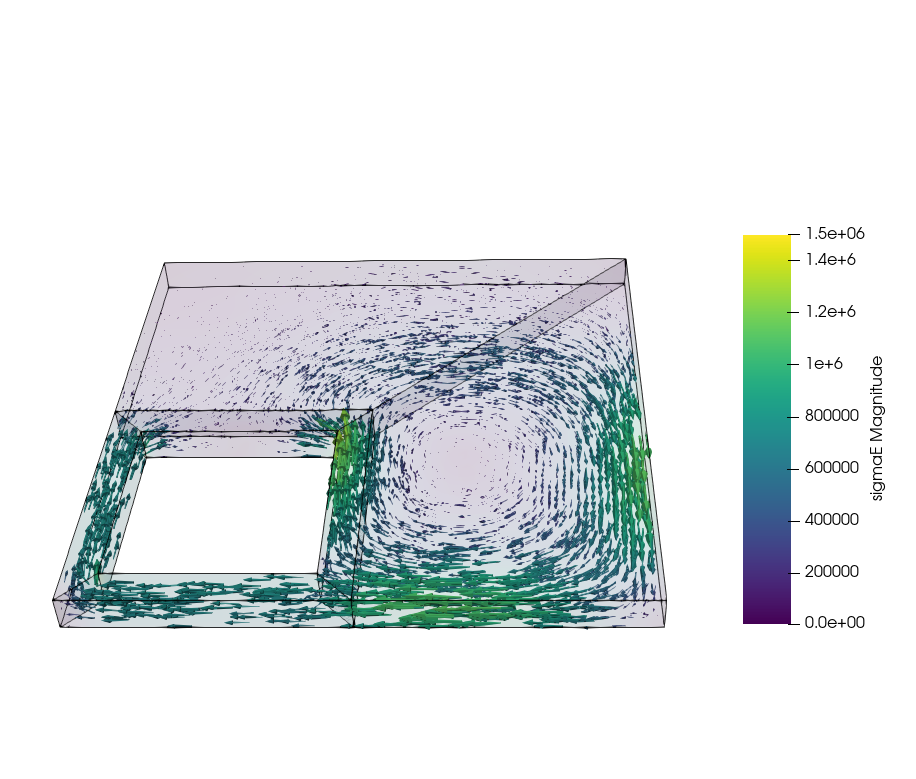}};
             \node [anchor=west] at (0,1.7) {\phantom{$\lvert\vecs{J}_\textrm{c}\rvert$ (\si{A\per\meter\squared})}};
        \end{tikzpicture}
        \\
        \begin{tikzpicture}
            \node at (0,0) {\includegraphics[height=2.5cm,trim=0cm 5cm 7cm 9cm, clip]{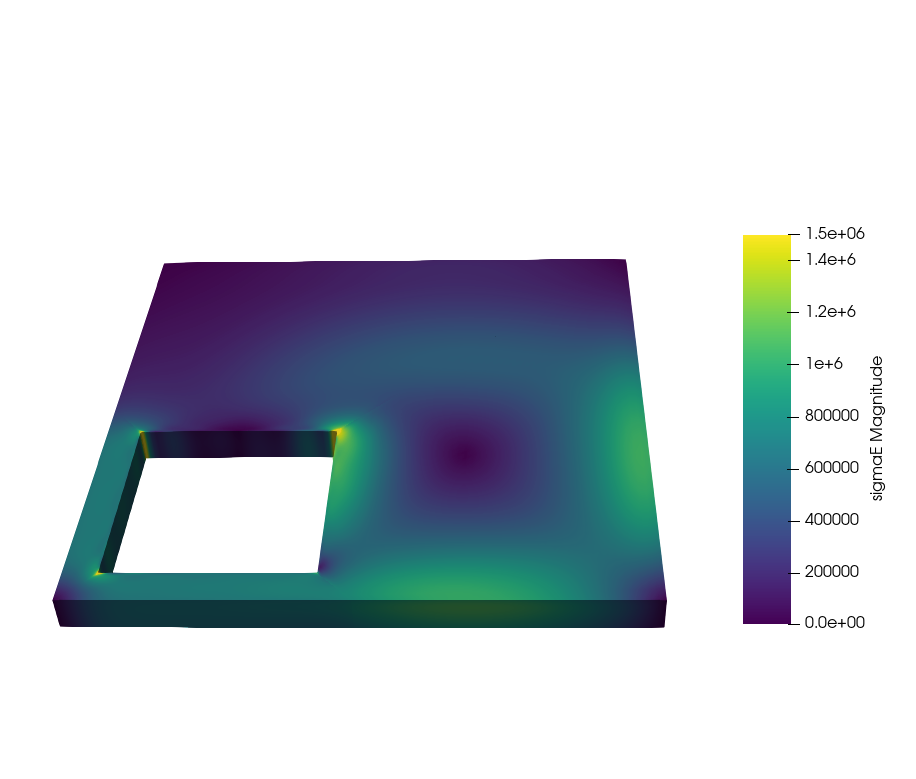}};
             \node [anchor=west] at (0,1.7) {\phantom{$\lvert\vecs{J}_\textrm{c}\rvert$ (\si{A\per\meter\squared})}};
        \end{tikzpicture}
        \subcaption{$\vecs{T}-\Omega$}
    \end{subfigure}
    \hfill
    \begin{subfigure}[c]{0.44\columnwidth}
        \center
        \begin{tikzpicture}
            \node at (0,0) {\includegraphics[height=2.5cm,trim=0cm 5cm 5cm 9cm, clip]{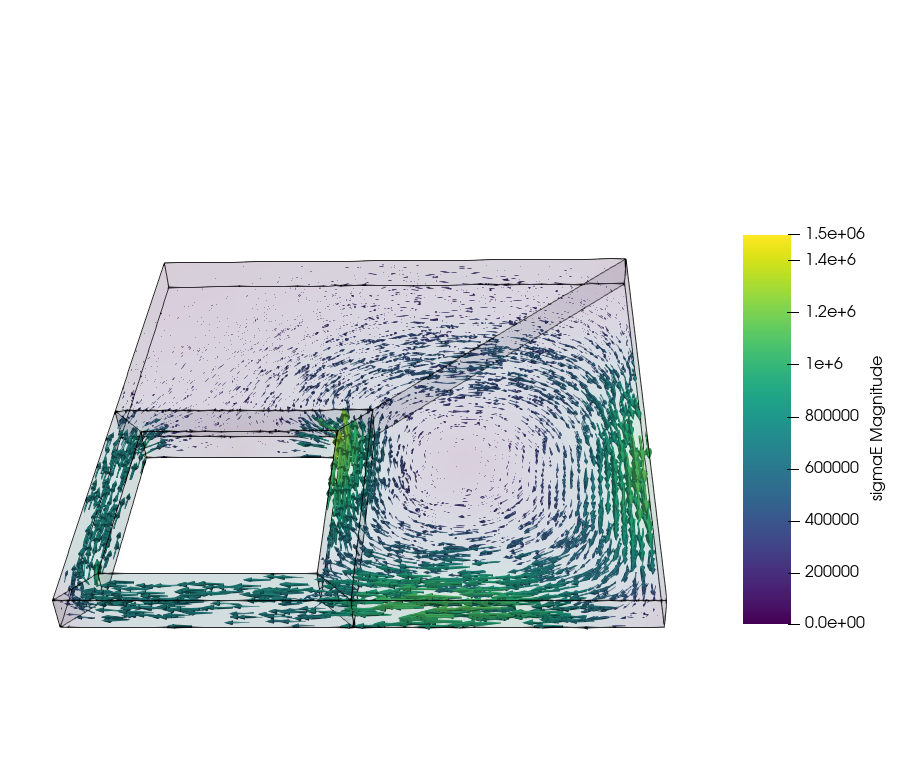}};
            \foreach \i in {0,0.3,...,1.51}{
                \draw (2.45,4/3*1.2*\i-1.15)--++(5pt,0) node [pos=1,anchor=west] {\small$\pgfmathprintnumber[fixed,precision=1]{\i}\cdot 10^6$};
             }
             \node [anchor=west] at (2.2,1.7) {$\lvert\vecs{J}_\textrm{c}\rvert$ (\si{A\per\meter\squared})};
        \end{tikzpicture}
        \\
        \begin{tikzpicture}
            \node at (0,0) {\includegraphics[height=2.5cm,trim=0cm 5cm 5cm 9cm, clip]{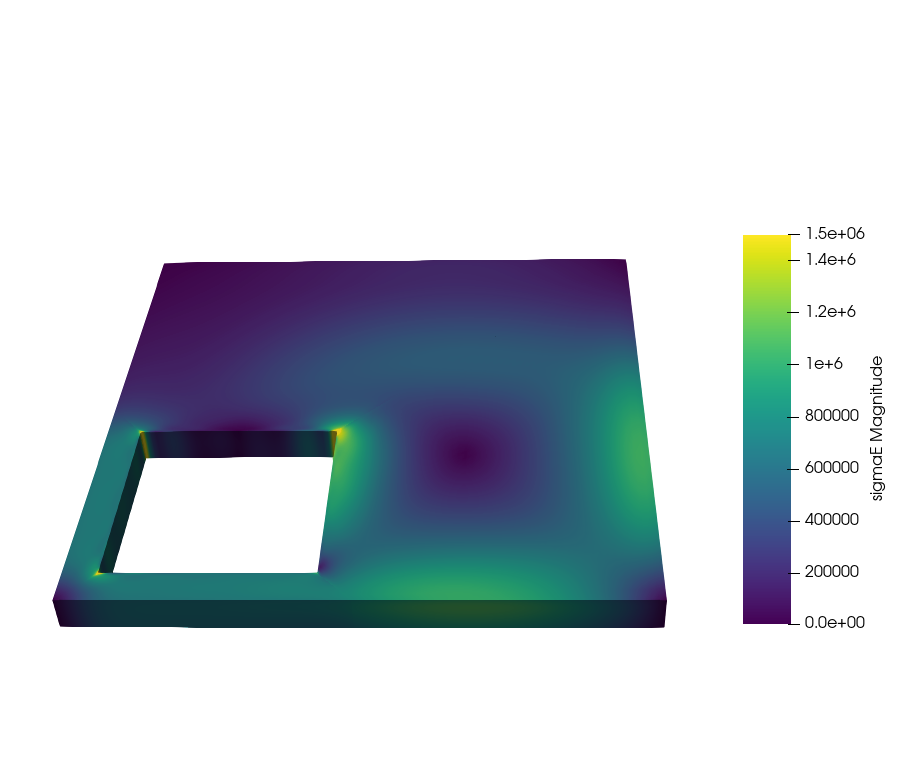}};
            \foreach \i in {0,0.3,...,1.51}{
                \draw (2.45,4/3*1.2*\i-1.15)--++(5pt,0) node [pos=1,anchor=west] {\small$\pgfmathprintnumber[fixed,precision=1]{\i}\cdot 10^6$};
             }
             \node [anchor=west] at (2.2,1.7) {$\lvert\vecs{J}_\textrm{c}\rvert$ (\si{A\per\meter\squared})};
        \end{tikzpicture}
        \subcaption{$\vecs{A}-\varphi$}
    \end{subfigure}
    \caption{Eddy currents $\vecs{J}_\mathrm{c}$ in the conducting plate at $t=\SI{0}{s}$ computed using the formulations in \cref{subsec:weak_IGA} with basis functions of degree $p=4$ . \textit{Top:} Direction of the currents $\vecs{J}_\textrm{c}$. \textit{Bottom:} Magnitude $\lvert\vecs{J}_\textrm{c}\rvert$ of the eddy currents. In Figure~\ref{fig:team7_J_Hphi} the evaluation line of Figure~\ref{fig:team7_Jy_line} is shown in red.} 
    \label{fig:team7_J}
\end{figure}

\begin{figure}
    \centering
    \begin{subfigure}[c]{0.49\columnwidth}
        \center
        \begin{tikzpicture}

              \begin{axis}[
                unbounded coords=jump, 
                ylabel={$J_y$ (\si{A\per\meter\squared})},
                xlabel={$x$  (\si{\meter})}, 
                width =1\textwidth, 
ylabel near ticks, 
                xlabel near ticks,
legend pos = south west,
                legend style={font=\small},
]
                        \addplot  [black, only marks
                        ] table [x index=0,  y expr=sign (\thisrowno{1})*sqrt((\thisrowno{2} )^2+ (\thisrowno{1})^2) , col sep=comma] {images/data/team7_solution_fujiwara_x_Jy_0_90.csv};
                        \addplot  [colorA, thick
                        ] table [x index=0, y index=1, col sep=comma] {images/data/team7_AV_A_evaluation_line_nsub_6_deg_3_refined.csv};
                        \addplot  [colorB,  thick,
                        ] table [x index=0, y index=1, col sep=comma] {images/data/team7_Hphi_evaluation_line_nsub_6_deg_3_refined.csv};
                        \addplot  [colorD, thick, dashed,
                        ] table [x index=0, y index=1, col sep=comma] {images/data/team7_TOmega_evaluation_line_nsub_6_deg_3_refined.csv};
                        \legend{
                        {measured \cite{Fujiwara_1990aa}},
                        {$\vecs{A}-\varphi$},
                        {$\vecs{H}-\phi$,}, 
                        {$\vecs{T}-\Omega$}, 
}
                \end{axis}
\end{tikzpicture}

         \subcaption{$p=3$}\label{fig:team7_Jy_line_p3}
    \end{subfigure}
    \begin{subfigure}[c]{0.49\columnwidth}
        \center
        \begin{tikzpicture}

              \begin{axis}[
                unbounded coords=jump, 
                ylabel={$J_y$ (\si{A\per\meter\squared})},
                xlabel={$x$  (\si{\meter})}, 
                width =1\textwidth, 
ylabel near ticks, 
                xlabel near ticks,
legend pos = south west,
                legend style={font=\small},
]
                        \addplot  [black, only marks
                        ] table [x index=0,  y expr=sign (\thisrowno{1})*sqrt((\thisrowno{2} )^2+ (\thisrowno{1})^2) , col sep=comma] {images/data/team7_solution_fujiwara_x_Jy_0_90.csv};
                        \addplot  [colorA, thick
                        ] table [x index=0, y index=1, col sep=comma] {images/data/team7_AV_A_evaluation_line_nsub_6_deg_4_refined.csv};
                        \addplot  [colorB,  thick,
                        ] table [x index=0, y index=1, col sep=comma] {images/data/team7_Hphi_evaluation_line_nsub_6_deg_4_refined.csv};
                        \addplot  [colorD, thick, dashed,
                        ] table [x index=0, y index=1, col sep=comma] {images/data/team7_TOmega_evaluation_line_nsub_6_deg_4_refined.csv};
                        \legend{
                        {measured \cite{Fujiwara_1990aa}},
                        {$\vecs{A}-\varphi$},
                        {$\vecs{H}-\phi$,}, 
                        {$\vecs{T}-\Omega$}, 
}
                \end{axis}
\end{tikzpicture}

         \subcaption{$p=4$}\label{fig:team7_Jy_line_p4}
    \end{subfigure}
    \caption{Comparison of the $y$-component of the eddy currents evaluated on the line $y=\SI{72}{mm}$, $z=\SI{19}{mm}$ on the conductor computed with different methods including the well-known $\vecs{A}-\phi$-formulation as a reference and the measured values of \cite{Fujiwara_1990aa}, for degree $p=3$ (\cref{fig:team7_Jy_line_p3}) and $p=4$ (\cref{fig:team7_Jy_line_p4}).}\label{fig:team7_Jy_line}
\end{figure}
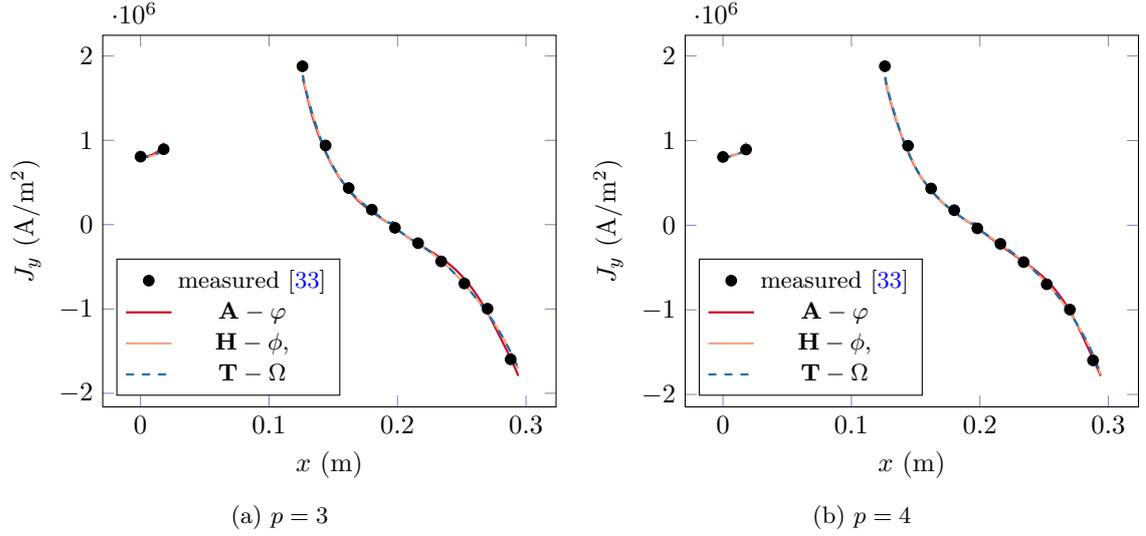

\begin{figure}
    \center
\begin{tikzpicture}
            \node at (0,0) {\includegraphics[width=0.5\columnwidth,trim=0cm 5cm 7cm 9cm, clip]{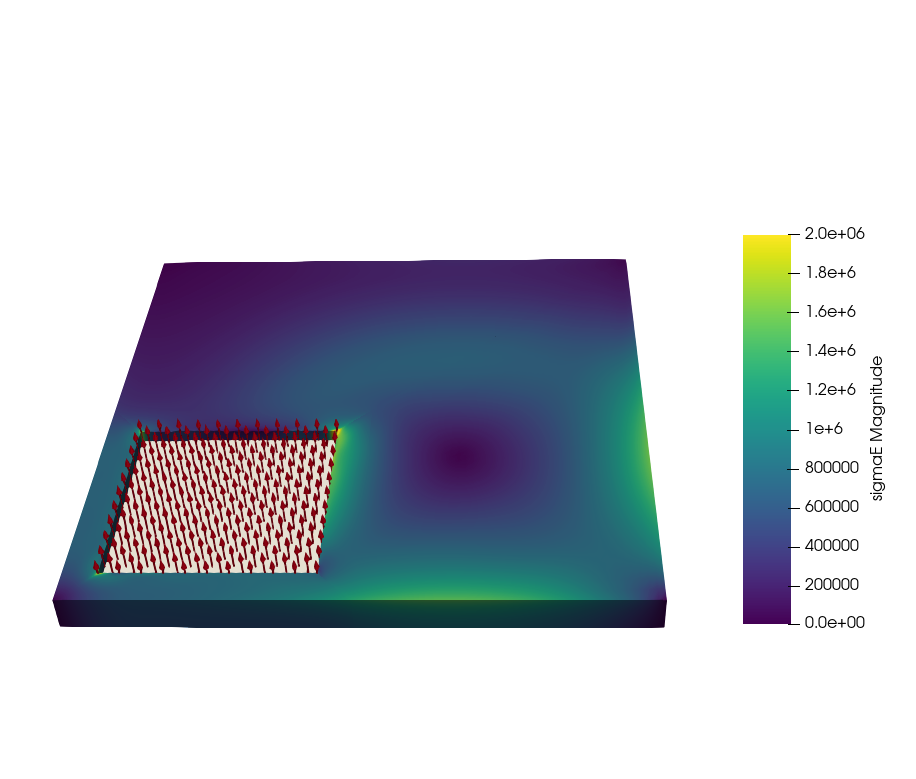}};
\end{tikzpicture}
    \caption{Support of the computed basis function of $W_h(\comdom_{\mathrm{i}})$, spanned by the red arrows.}
    \label{fig:team7_harmonic_fields}
\end{figure}

The eddy currents $\vecs{J}_\textrm{c} = \sigma \vecs{E}$ that are excited in the conductor computed using the $\vecs{H}-\phi$, $\vecs{T}-\Omega$ and $\vecs{A}-\varphi$ formulations are shown in Figure~\ref{fig:team7_J}. The whole domain is discretized by $39$ patches and a total of $8892$ elements. The conductor consists of $6$ patches and $1080$ elements, and the total number of unknowns is $26561$ for the $\vecs{H}-\phi$, $\vecs{T}-\Omega$ formulations, and $61649$ for the $\vecs{A}-\varphi$ formulation when discretizing with basis functions of degree $p=3$ and $37738$ for the $\vecs{H}-\phi$, $\vecs{T}-\Omega$ formulations, and $87976$ for the $\vecs{A}-\varphi$ formulation when discretizing with basis functions of degree $p=4$. Since there is a single hole, there is only one basis function in $W_h(\comdom_{\mathrm{i}})$. 
The number $\Nharm$ of basis functions of $S_p^1(\comdom_{\mathrm{i}})$ which are linearly combined to obtain this function, and the number $\Ngrad$ of basis functions of $S_p^0(\comdom_{\mathrm{i}})$, are $\Nharm=81$ and $\Ngrad=20218$ for the case of $p=3$ and $\Nharm=100$ and $\Ngrad=28262$ for the case of $p=4$.
To quantitatively compare the results of the different formulations we evaluate the $y$-component $J_y$ of the eddy currents $\vecs{J}_\mathrm{c}$ along the evaluation line used in \cite{Ledger_2010aa,Fujiwara_1990aa} on the conductor, i.e., $y=\SI{72}{mm}$, $z=\SI{19}{mm}$ shown in Figure~\ref{fig:team7_Jy_line}. Here, the eddy currents are computed as a complex quantity. Following \cite{Fujiwara_1990aa}, we consider for complex quantities
\begin{align*}
 \chi = \operatorname{sign}(\chi_r) \sqrt{\chi_r^2+\chi_i^2},
\end{align*}
where $\chi_r$ and $\chi_i$ are respectively the real and imaginary parts of $\chi$. The results obtained using the different formulations are compared to the measured results given in \cite{Fujiwara_1990aa}.
We can see, that the results are in good agreement with the measured eddy currents on the conductor. To show qualitatively how the cohomology generator influences the fill in of the system we additionally show how the vector field generated by the basis function in $W_h(\comdom_{\mathrm{i}})$. The field has indeed very local support inside the conductor hole as shown in Figure~\ref{fig:team7_harmonic_fields}.
\subsection{Robustness with respect to cohomology group dimension}
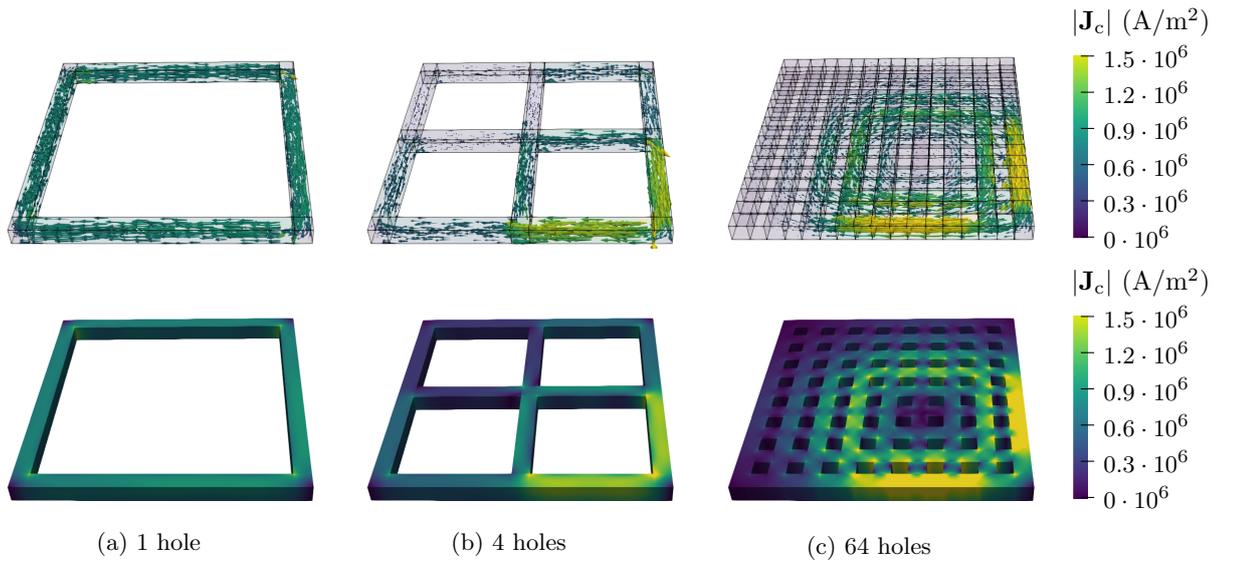
\begin{figure}
    \center
    \begin{subfigure}[c]{0.3\columnwidth}
        \center
        \begin{tikzpicture}
            \node at (0,0) {\includegraphics[height=2.5cm,trim=0cm 5cm 7cm 9cm, clip]{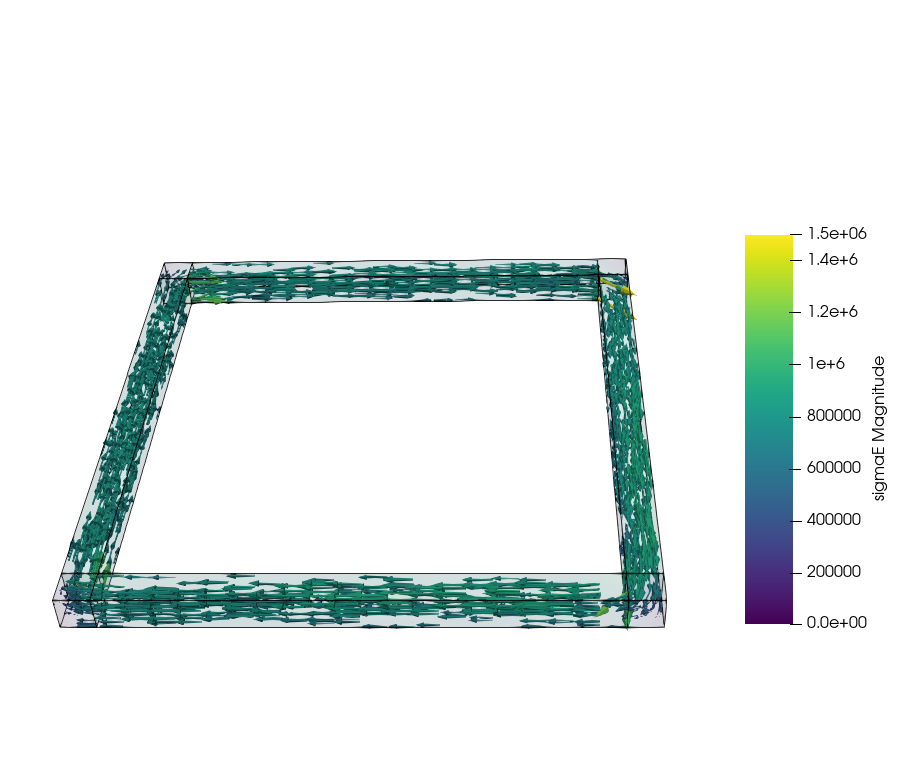}};
            \node [anchor=west] at (0,1.7) {\phantom{$\lvert\vecs{J}_\textrm{c}\rvert$ (\si{A\per\meter\squared})}};
        \end{tikzpicture}
        \\
        \begin{tikzpicture}
            \node at (0,0) {\includegraphics[height=2.5cm,trim=0cm 5cm 7cm 9cm, clip]{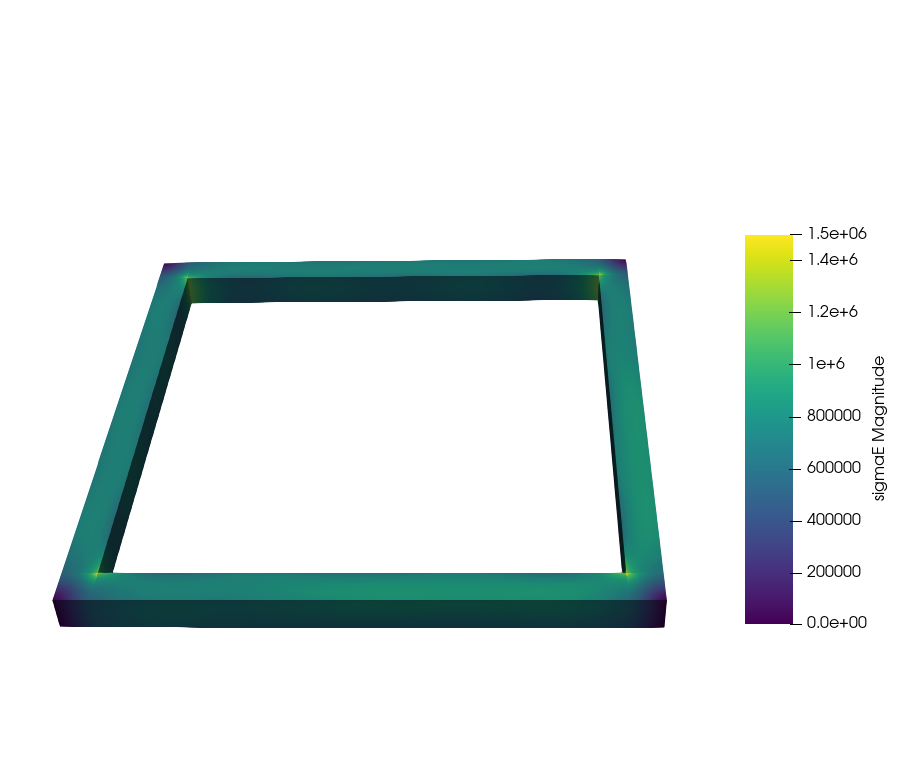}};
            \node [anchor=west] at (0,1.7) {\phantom{$\lvert\vecs{J}_\textrm{c}\rvert$ (\si{A\per\meter\squared})}};
        \end{tikzpicture}
        \subcaption{$1$ hole}\label{fig:manyHoles_1}
    \end{subfigure}
    \begin{subfigure}[c]{0.3\columnwidth}
        \center
        \begin{tikzpicture}
            \node at (0,0) {\includegraphics[height=2.5cm,trim=0cm 5cm 7cm 9cm, clip]{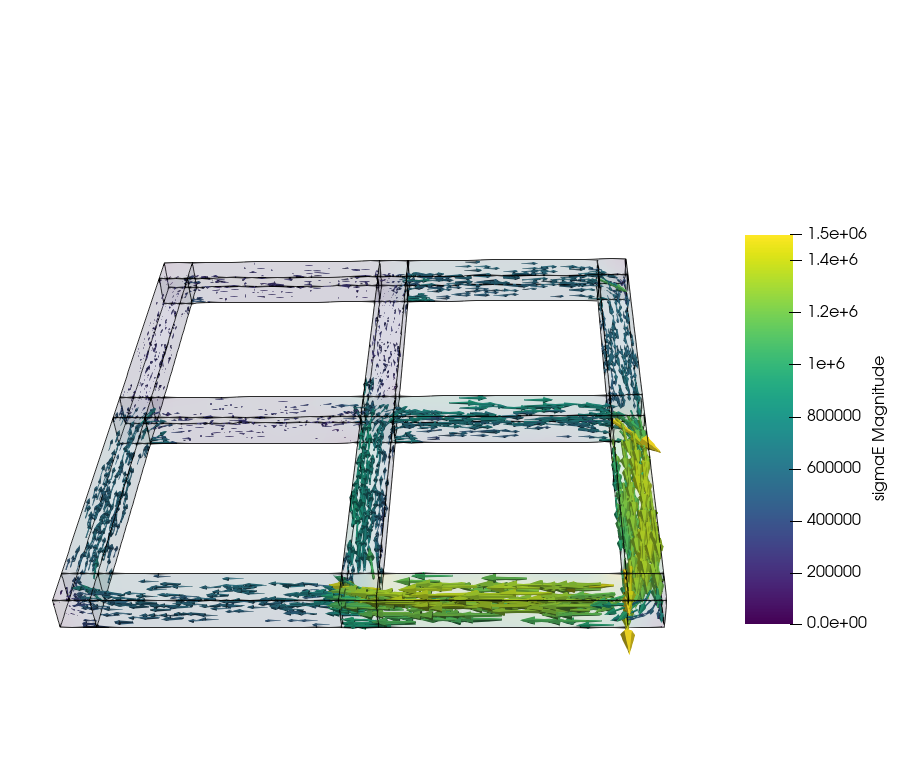}};
             \node [anchor=west] at (0,1.7) {\phantom{$\lvert\vecs{J}_\textrm{c}\rvert$ (\si{A\per\meter\squared})}};
        \end{tikzpicture}
        \\
        \begin{tikzpicture}
            \node at (0,0) {\includegraphics[height=2.5cm,trim=0cm 5cm 7cm 9cm, clip]{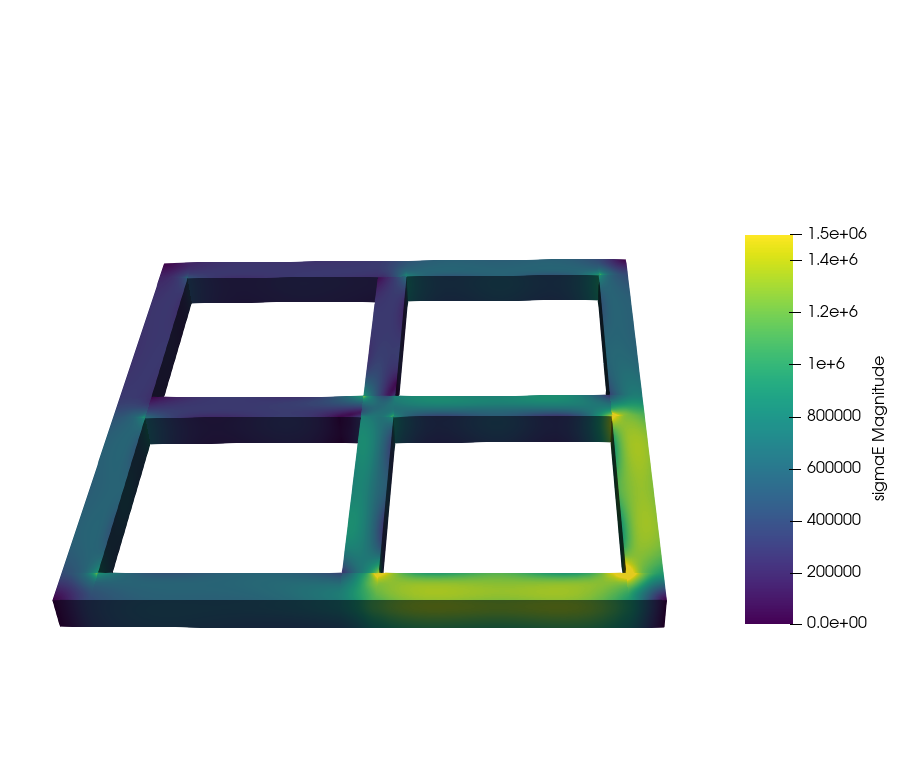}};
             \node [anchor=west] at (0,1.7) {\phantom{$\lvert\vecs{J}_\textrm{c}\rvert$ (\si{A\per\meter\squared})}};
        \end{tikzpicture}
        \subcaption{$4$ holes}\label{fig:manyHoles_4}
    \end{subfigure}
    \begin{subfigure}[c]{0.3\columnwidth}
        \center
        \begin{tikzpicture}
            \node at (0,0) {\includegraphics[height=2.5cm,trim=0cm 5cm 5cm 9cm, clip]{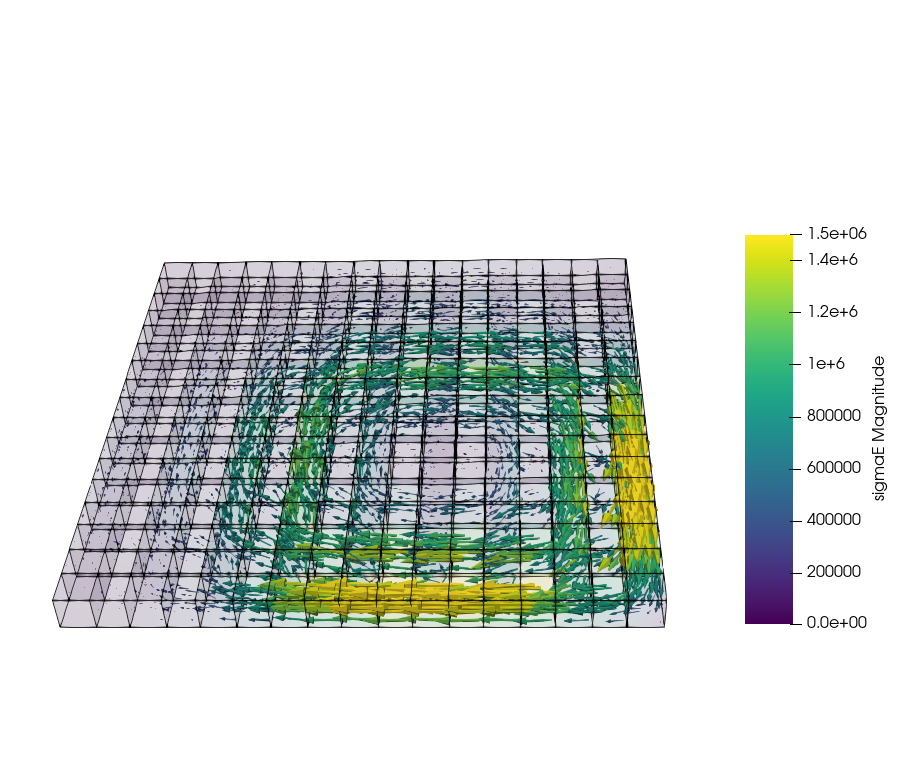}};
            \foreach \i in {0,0.3,...,1.51}{
                \draw (2.45,4/3*1.2*\i-1.15)--++(5pt,0) node [pos=1,anchor=west] {\small$\pgfmathprintnumber[fixed,precision=1]{\i}\cdot 10^6$};
             }
             \node [anchor=west] at (2.2,1.7) {$\lvert\vecs{J}_\textrm{c}\rvert$ (\si{A\per\meter\squared})};
        \end{tikzpicture}
        \\
        \begin{tikzpicture}
            \node at (0,0) {\includegraphics[height=2.5cm,trim=0cm 5cm 5cm 9cm, clip]{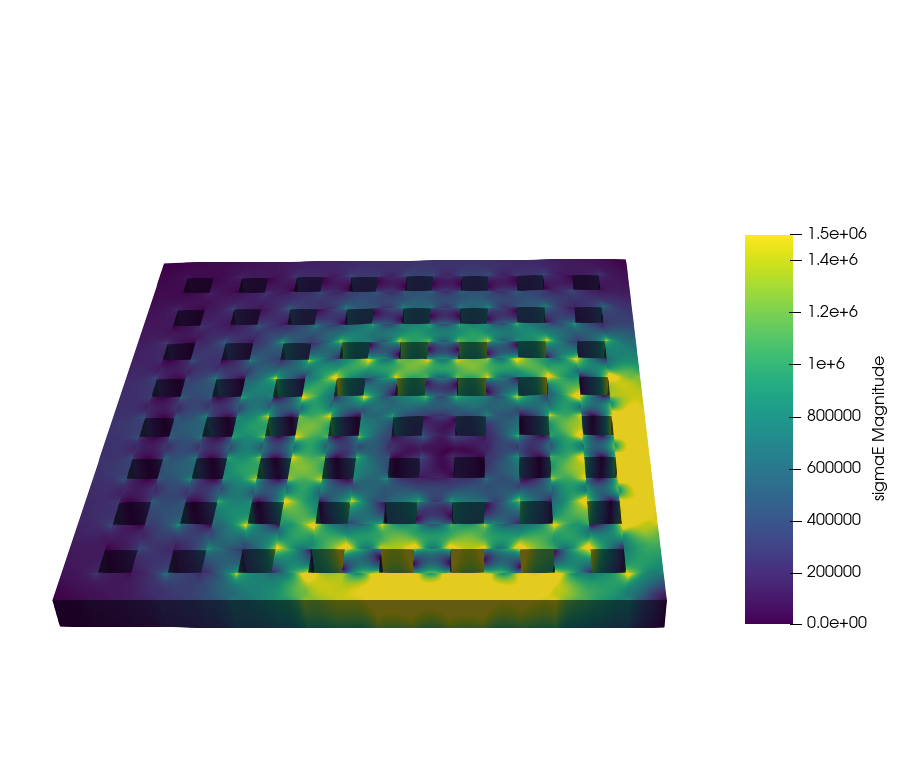}};
            \foreach \i in {0,0.3,...,1.51}{
                \draw (2.45,4/3*1.2*\i-1.15)--++(5pt,0) node [pos=1,anchor=west] {\small$\pgfmathprintnumber[fixed,precision=1]{\i}\cdot 10^6$};
             }
             \node [anchor=west] at (2.2,1.7) {$\lvert\vecs{J}_\textrm{c}\rvert$ (\si{A\per\meter\squared})};
        \end{tikzpicture}
        \subcaption{$64$ holes}\label{fig:manyHoles_64}
    \end{subfigure}
    \caption{Eddy currents $\vecs{J}_\mathrm{c}$ in the conducting plate with varying number of holes at $t=\SI{0}{\second}$. Computed using the $\vecs{H}-\phi$-formulation with degree $p=3$. \textit{Top:} Direction of the currents $\vecs{J}_\textrm{c}$. \textit{Bottom:} Magnitude $\lvert\vecs{J}_\textrm{c}\rvert$ of the eddy currents.} 
    \label{fig:manyHoles_J}
\end{figure}

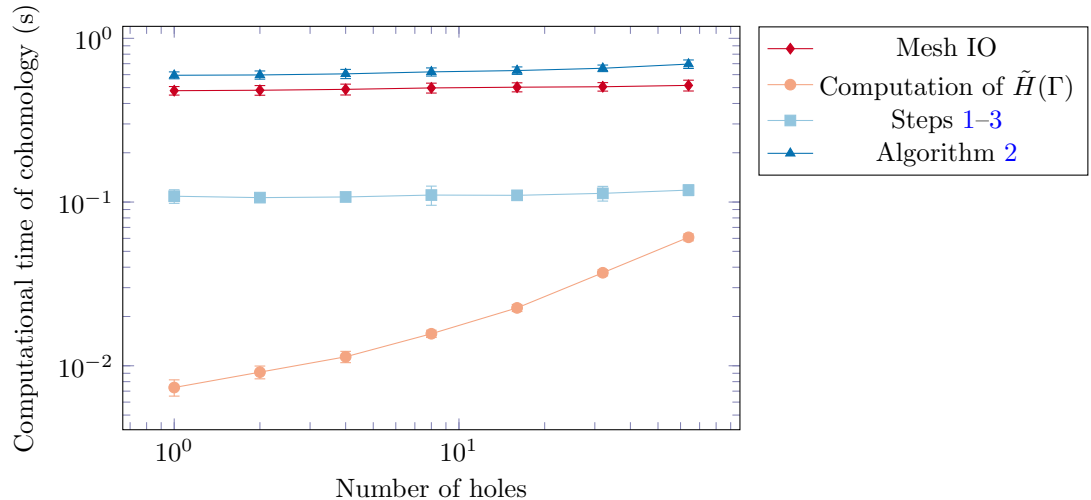
\begin{figure}
    \centering
    \begin{tikzpicture}
                \begin{loglogaxis}[
                ylabel={Computational time of cohomology (s)},
                xlabel={Number of holes}, 
                width = 0.63\columnwidth, 
                height=0.45\columnwidth, 
                ylabel near ticks, 
                xlabel near ticks,
                legend pos=outer north east,
]

\addplot+[colorA,mark=diamond*, mark options={fill=colorA},error bars/.cd,y dir=both, y explicit]
                        table
                        [col sep=comma, 
                        x index=0, 
                        y expr=\thisrowno{1}, 
                        y error expr=\thisrowno{7}] {images/data/topoprotimings100samples_holes_mean_standardDeviation.csv};
\addplot+[colorB,mark=*, mark options={fill=colorB},error bars/.cd,y dir=both, y explicit]
                        table
                        [col sep=comma, 
                        x index=0, 
                        y expr=\thisrowno{2}, 
                        y error expr=\thisrowno{8}] {images/data/topoprotimings100samples_holes_mean_standardDeviation.csv};
                        
\addplot+[colorC,mark=square*, mark options={fill=colorC},error bars/.cd,y dir=both, y explicit]
                        table
                        [col sep=comma, 
                        x index=0, 
                        y expr=\thisrowno{3}, 
                        y error expr=\thisrowno{9}] {images/data/topoprotimings100samples_holes_mean_standardDeviation.csv};
                        
\addplot+[colorD,mark=triangle*, mark options={fill=colorD}, error bars/.cd,y dir=both, y explicit]
                        table
                        [col sep=comma, 
                        x index=0, 
                        y expr=\thisrowno{5}, 
                        y error expr=\thisrowno{11}] {images/data/topoprotimings100samples_holes_mean_standardDeviation.csv};

                      \legend{Mesh IO, Computation of $\tilde{H}(\Gamma)$, Steps \ref{algo:step1}--\ref{algo:step3}, Algorithm \ref{algo:generators}}
                \end{loglogaxis}
\end{tikzpicture}

     \caption{Computational time and standard deviation of the cohomology computation depending on the number of holes in the conductor.}\label{fig:topopro_time}
\end{figure}

To investigate the influence of the topology on the computational time of the cohomology basis functions, we construct an example similar to the TEAM benchmark problem 7. The excitation is the same as in the TEAM problem 7 (see \cite{Fujiwara_1990aa}), however, the conductor has a variable number of holes, where the width of the conducting parts is always $\SI{18}{mm}$. The eddy currents for a few configurations are shown in Figure~\ref{fig:manyHoles_J}.
We finally use an additional manufactured example setup to show that for the same mesh the computational cost for the cohomology generators is very robust with respect to the number of holes in the conductor. To investigate this, we measure the time for the computation of the cohomology generators for a variable number of holes in the conductor. The discretization, i.e. number of patches/elements stays the same throughout the computation and only the material distribution is changed. In each case, the geometry consists of $1083$ patches and a total of $1083$ elements and basis functions of degree $p=3$ are used. For each setting, the computation is repeated $100$ times and the mean computational time and standard deviation are evaluated and shown in Figure~\ref{fig:topopro_time}.
{The results show in logarithmic scale that the only part of the computation which does not stay constant is the computation of generators for  $\tilde{H}(\Gamma_C)$, that is, the cohomology generators of the interface,
which is a preliminary step to \cref{algo:generators}. This is due to the fact that, even though the mesh for the whole domain is unchanged (as shown by the fact that the cost for mesh input--output between GeoPDEs and Topoprocessor stays constant) the interface mesh grows in size. This is clearly visible going from left to right in Figure~\ref{fig:manyHoles_J}. Thus, the graph associated to the mesh discretization of interface $\Gamma$, over which STT algorithms behave linearly with respect to number of edges, grows. The last curve in Figure~\ref{fig:topopro_time} shows the total cost of the algorithm, summing up all contributions, and highlights the fact that the leading cost is by far building the cell complex for the cohomology computations. Its overhead can nevertheless be minimized by a streamlined implementation of the whole formulation within a single toolbox, which goes beyond the goals of the present article.}
 
\section{Conclusions}\label{sec:conclusions}
Based on the spline-version of the de Rham complexes approximation~\cite{Beirao-da-Veiga_2014aa}, we show that the algorithms of~\cite{Dlotko_2017aa} can be extended to isogeometric analysis with high order splines. We exploit the the isomorphisms between the spline spaces and low order finite elements defined on the hexahedral control mesh, to apply existing algorithms from finite elements to this auxiliary control mesh. This enables automated cohomology computations for spline-based discretizations of the eddy current problem in magnetic formulations. The numerical results confirm the convergence of the proposed method and the optimal scaling of the algorithms with respect to the dimension of the cohomology group of the underlying domain, which is very much desired when complicated devices with hundreds or thousands of holes have to be simulated (e.g. in electric motors or magnetic confinement of plasma in tokamaks). 
\section*{Acknowledgements}{This work is supported by the Graduate School CE within the Centre for Computational Engineering at Technische Universität Darmstadt, by the Swiss National Science Foundation via the project HOGAEMS n.200021\_188589, by the German Research Foundation (DFG) via the project SCHO 1562/6-1 and the joint DFG/FWF Collaborative Research Centre CREATOR (CRC – TRR361/F90) at TU Darmstadt, TU Graz and JKU Linz.

\end{document}